\documentclass[preprint,12pt,times]{elsarticle}

\usepackage{amssymb}
\usepackage{todonotes}
\usepackage[acronym]{glossaries}
\usepackage{comment}
\usepackage{cleveref}
\usepackage{tikz}
\usepackage{amsmath}
\usepackage{amsthm}
\usepackage{mathtools}
\usepackage{siunitx}
\usetikzlibrary{shapes,arrows}
\usetikzlibrary{positioning}
\usetikzlibrary{fit}
\usetikzlibrary{decorations.pathreplacing}
\usepackage{bm}

\newcounter{dummy} \numberwithin{dummy}{section}

\newtheorem{theorem}[dummy]{Theorem}

\journal{Journal of Computational Physics}

\makeglossaries
\newacronym{UQ}{UQ}{Uncertainty Quantification}
\newacronym{MC}{MC}{Monte Carlo}
\newacronym{MCMC}{MCMC}{Markov chain Monte Carlo}
\newacronym{MLMCMC}{MLMCMC}{Multilevel Markov chain Monte Carlo}
\newacronym{MLMC}{MLMC}{Multilevel Monte Carlo}
\newacronym{HPC}{HPC}{High Performance Computing}
\newacronym{DUNE}{DUNE}{Distributed and Unified Numerics Environment}
\newacronym{GMsFEM}{GMsFEM}{Generalized Multiscale Finite Element Method}
\newacronym{GEVP}{GEVP}{generalized eigenvalue problem}
\newacronym{MS-GFEM}{MS-GFEM}{multiscale-spectral generalized finite element method}
\newacronym{ATI}{ATI}{Aerospace Technology Institute}
\newacronym{GenEO}{GenEO}{Generalized Eigenproblems in the Overlaps}
\newacronym{A-h-GenEO}{A-h-GenEO}{A-harmonic GenEO}
\newacronym{FE}{FE}{finite element}
\newacronym{DoF}{DoF}{degree of Freedom}
\newacronym{PDE}{PDE}{partial differential equation}
\newacronym{LMOR}{LMOR}{localized model order reduction}
\newacronym{UMFPACK}{UMFPACK}{Unsymmetric MultiFrontal PACKage}
\newacronym{CG}{CG}{Conjugate Gradient}
\newacronym{GMRES}{GMRES}{Generalized Minimal Residual}
\newacronym{AMG}{AMG}{Algebraic Multigrid}
\newacronym{DD}{DD}{Domain Decomposition}
\newacronym{MAP}{MAP}{Maximum a Posteriori Probability}
\newacronym{o}{o}{overlap size}
\newacronym{ovs}{ovs}{oversampling size}
\newacronym{REV}{REV}{representative elementary volume}
\newacronym{PoU}{PoU}{partition of unity}

\begin{document}

\begin{frontmatter}



\title{Scalable multiscale-spectral GFEM with an application to composite aero-structures}


\author[inst1]{Jean Bénézech\corref{cor1}}
\cortext[cor1]{Corresponding author}
\ead{jb3285@bath.ac.uk}

\author[inst2]{Linus Seelinger}

\author[inst2]{Peter Bastian}

\author[inst1]{Richard Butler}

\author[inst3]{Timothy Dodwell}

\author[inst7]{Chupeng Ma}

\author[inst2,inst5]{Robert Scheichl}

\affiliation[inst1]{organization={Centre for Integrated Materials, Processes and Structures},
            addressline={University of Bath}, 
            city={Bath},
            postcode={BA2 7AY}, 
            country={UK}}

\affiliation[inst2]{organization={Institute for Applied Mathematics and Interdisciplinary Center for Scientific Computing, Heidelberg University},
            city={Heidelberg},
            country={Germany}}

\affiliation[inst3]{organization={Institute of Data Science and AI, University of Exeter},
            city={Exeter},
            country={UK \ \& \ digiLab, Exeter, UK}}


\affiliation[inst5]{organization={Department of Mathematical Sciences, University of Bath},
            city={Bath},
            country={UK}}
            

\affiliation[inst7]{organization={Institute of Scientific Research, Great Bay University},
            city={Dongguan},
            country={China}}

\begin{abstract}

In this paper, the first large-scale application of multiscale-spectral generalized finite element methods (MS-GFEM) to composite aero-structures is presented. 
The crucial novelty lies in the introduction of A-harmonicity in the local approximation spaces, which in contrast to [Babuska, Lipton, \emph{Multiscale Model.~Simul.}~\textbf{9}, 2011] is enforced more efficiently via a constraint in the local eigenproblems.   
This significant modification leads to excellent approximation properties, which turn out to be essential to capture accurately material strains and stresses with a low dimensional approximation space, hence maximising model order reduction.
The implementation of the framework in the \gls{DUNE} software package, as well as a detailed description of all components of the method are presented and exemplified on a composite laminated beam under compressive loading.
The excellent parallel scalability of the method, as well as its superior performance compared to the related, previously introduced GenEO method are demonstrated on two realistic application cases, including a C-shaped wing spar with complex geometry.
Further, by allowing low-cost approximate solves for closely related models or geometries this efficient, novel technology provides the basis for future applications in optimisation or uncertainty quantification on challenging problems in composite aero-structures.
\end{abstract}



\begin{keyword}
Multiscale-spectral Generalized Finite Element Methods (MS-GFEM) \sep GenEO coarse space \sep A-harmonic subspace \sep large-scale composite structure \sep parallel scalability
\end{keyword}

\end{frontmatter}


\section{Introduction}

Ill-conditioned and multiscale partial differential equations (PDEs) arise in many fields, where the  computation of a resolved, fine-scale solution or a robust low-dimensional approximation can be challenging.
Due to the interaction of their mesoscopic structure (ply level; sub-millimeter scale) and their geometric macroscopic features (structural level; meter scale), composite aero-structures are naturally, inherently multiscale. To model the behavior of a large scale composite structure, the numerical model needs to accurately represent the meso-scale configuration of the material as well as the macro-scale geometry of the part. Full mesoscopic descriptions of large components naturally lead to models with huge numbers of degrees of freedom. This makes the computations prohibitively expensive, particularly in contexts where many evaluations are required, e.g. optimisation or uncertainty quantification. As a result, composite aero-structures provide an ideal test bed for the new multiscale method proposed in this paper, which allows for the interaction of fine and coarse scale behaviour to be captured without becoming excessive in cost.

\subsection{High Performance Solvers for Composite Applications}

Due to strongly varying parameters across the simulation domain, elasticity problems arising in composite materials lead to Finite Element (FE) matrices that are extremely ill-conditioned \cite{reinarz2018dune,butler2020high}. This poor conditioning is due to the contrast in stiffness between the carbon fibres and the surrounding resin matrix, as well as the complex anisotropy arising from the inclusion of long directional fibres. As a result, composite laminates have both low-energy modes, whereby the stiff fibres act as rigid body inclusions and the complement resin deforms easily; and also very high-energy modes, in deformation regimes with stretched, stiff fibres. This contrast between low and high energy modes is at the heart of the ill-conditioning (high condition number) in composite applications. 

Whilst sparse direct solvers, like UMFPACK 
\cite{UMFPACK} or the ones provided in the \emph{Abaqus} package \cite{systemes2007abaqus}, can reliably solve such systems, they are inherently limited in their scalability. This immediately restricts the physical scale of composites that can be simulated. Iterative solvers such as \gls{CG}
or GMRES 
\cite{Saad2003iterative}
in turn promise massive parallel scalability for modern \gls{HPC} systems, but their efficiency (i.e. number of iterations) strongly depends on the condition number of the matrix.

In order to render such iterative solvers robust, preconditioners are essential. Whilst \gls{AMG} preconditioners are in general a promising choice regarding robustness and scalability, tests with two AMG implementations, in dune-istl \cite{blatt2006iterative} and BoomerAMG \cite{boomerAMG}, have demonstrated poor performance in composite applications \cite{reinarz2018dune,butler2020high}. The reason is that, without a problem-specific local aggregation strategy, coarse grids in AMG do not capture the low-energy modes in composites structures, which motivates the search for alternative strategies that produce fast approximate coarse representations at a larger scale.
The \gls{GenEO} method \cite{spillane2014abstract} provides such a coarsening, leading to a robust, theoretical bound on the condition number of the preconditioned system for a two-level Schwarz \gls{DD} method \cite{toselli2004domain}. A problem-specific coarse space is computed from some tailored, local eigenproblems on overlapping subdomains. Robust scale-up to several thousands of processor cores for composites applications has been shown in \cite{butler2020high, geneo-hpc}. While the iterative solver now only needs few iterations, a considerable cost is expended in solving the independent local eigenproblems.

To avoid these tremendous computational costs of direct simulations of multiscale problems, computational homogenization methods \cite{guedes1990preprocessing,kouznetsova2004multi, nguyen2011multiscale} have been well developed and widely used in the engineering community. Moreover, in practical engineering applications, multiscale problems are typically solved multiple times with different source terms and possibly local changes in model parameters, such that the higher setup cost of such methods can be offset. Most of those methods, however, are based on scale separation hypotheses, and may fail for typical problems in realistic applications that do not exhibit such a scale separation.

\subsection{Multiscale Methods in Composite Analysis}
To efficiently solve multiscale problems without scale separation for repeated analysis required in an uncertainty quantification context, various multiscale model order reduction methods have been developed, such as the multiscale finite element method (MsFEM) \cite{hou1997multiscale}, the generalized finite element method (GMsFEM) \cite{efendiev2013generalized}, localized orthogonal decomposition (LOD) \cite{maalqvist2014localization}, flux norm homogenization \cite{berlyand2010flux}, and the generalized finite element method (GFEM) \cite{babuska2011optimal,babuvska2014machine,babuska2020multiscale}, to cite a few. 
Most of these methods were developed in the context of numerical multiscale methods and are based on representing the solution space on a coarse grid by a low dimensional space that is spanned by some pre-computed local basis functions that take into account the structural meso-scale information in the material parameters. 

The multiscale-spectral generalized finite element method (MS-GFEM), the focus of this work, was first proposed by Babuska and Lipton \cite{babuska2011optimal} for solving heterogeneous diffusion problems, but motivated also by problems in linear elasticity and in particular fibre-reinforced composites. The approach builds optimal local approximation spaces from eigenvectors of local eigenproblems posed on A-harmonic spaces defined for oversampling subdomains. Crucially, the global approximation error is fully controlled by the local approximation errors, which are rigorously proved to decay nearly exponentially -- a feature not shared by most other ad-hoc constructed numerical multiscale approaches.
For diffusion problems the implementation details 
were discussed in \cite{babuska2020multiscale} together with numerical results on a two-dimensional
toy example. 

In a recent paper \cite{ma2022novel}, Ma, Scheichl and Dodwell proposed new local eigenproblems involving the partition of unity to construct new optimal local approximation spaces for the MS-GFEM method, resulting in a GenEO-type coarse-space approximation. Significant advantages of the new local approximation spaces were demonstrated and sharper decay rates for the local approximation errors were proved. In \cite{ma2021error}, the MS-GFEM method was then also formulated and analysed for the first time in the discrete setting, as a non-iterative domain decomposition type method for solving linear systems resulting from FE discretizations of the fine-scale problem. Very similar local and global error estimates as in the continuous setting were derived. Furthermore, an efficient method to solve the (discrete) local eigenproblems was proposed, where the A-harmonic condition is directly incorporated into the eigenproblem. More recently, this approach was applied to other multiscale PDEs, such as Helmholtz \cite{ma2021wavenumber}, parabolic \cite{schleuss2022optimal} and singularly perturbed \cite{ma2022exponential} problems. 
Although numerical results of various simple, two-dimensional examples have demonstrated the efficiency of the MS-GFEM method, up to now, there is no study available in the literature on the application of the MS-GFEM method to realistic, large-scale, three-dimensional multiscale problems and on the implementation and performance of the method on massively parallel computers.

\subsection{Contributions of this paper}

This paper represents the first large-scale application of GenEO as a multiscale-spectral GFEM, providing good approximations to fine-scale solutions with a very low number of basis functions. 
The reformulation of GenEO as a GFEM method in local A-harmonic subspaces distinguishes our methodology from the one proposed by Babuska and Lipton \cite{babuvska2014machine}. Most notably, in contrast to \cite{babuvska2014machine} our method is inherently adaptive --- we can control the error a posteriori by simply setting a threshold on the eigenvalues to decide which eigenvectors need to be included into the local spaces; see Theorem~\ref{thm:error_bound}.
The paper constitutes an extension of the work proposed by Ma and co-authors \cite{ma2022novel,ma2021error} by generalizing the method to three dimensional elasticity problems, demonstrated with two real-world application cases.
The theoretical background of our formulation as well as its implementation are described in detail and its excellent performance and scalability are demonstrated. 
The resulting coarse space turns out to have significantly better approximation properties than GenEO in the elasticity problems considered here. In particular, local A-harmonicity is crucial for accurate strain approximation with significant practical improvements demonstrated in numerical experiments.

Our method provides efficient reduced order models for large-scale problems that exhibit strong dependence on local details. The efficiency arises from a decomposition of the global problem into independent sub-problems that can be treated fully in parallel. A scalability test is presented in this paper that demonstrates this high efficiency even for very large structures.
The accuracy of the approximation space is fully adjustable: a single threshold parameter on the local eigenvalues allows for an optimisation of the amount of model order of reduction. 
The approach does not rely on a scale separation hypothesis between material scales (meso- and macro-scale for the examples illustrated in the paper), the proposed multi-scale framework is particularly well suited for composite aero-structures.
This is demonstrated via an application on a realistic C-spar model -- a demonstrator application in the UK-EPSRC-funded CerTest project on composite structural design and certification (see the acknowledgements in Section \ref{sec:acknowledge}).

A key motivation for the approach presented in this paper is the development of an offline-online framework, where the costly local model order reductions are reused across multiple similar simulation runs, reducing the overall cost to a fraction, as suggested in \cite{efendiev2013generalized} but with significantly smaller coarse spaces. This promises to accelerate uncertainty quantification or optimisation tasks on challenging composites models and will be described in detail in a subsequent paper.

In addition, this paper reports on other technical improvements and recent developments of the \texttt{dune-composites} module within \gls{DUNE} \cite{DUNE}.
Support for GenEO is extended to unstructured DUNE grids, which is crucial for engineering applications such as composite parts with complex geometry (e.g., T-joint stiffeners) that cannot be discretized solely using a structured grid. Since the implementations of unstructured grids do not natively support overlapping subdomains in \gls{DUNE}, the overlapping FE matrices needed in GenEO are constructed directly from assembled ones on a non-overlapping grid partition.

\section{Problem formulation}

Let us start by formulating the anisotropic, linear elasticity equations for composite structures and their finite element discretization.
The composite structure is assumed to occupy a bounded and (for simplicity) polyhedral domain $\Omega \subset \mathbb R^3$ with boundary $\Gamma$ and unit, outward normal ${\bf n} \in \mathbb R^3$. At each point ${\bf x}\in \Omega$ we define a vector-valued displacement ${\bf u}({\bf x}): \Omega \rightarrow \mathbb R^3$ 
and denote by ${\bf f}({\bf x}): \Omega \rightarrow \mathbb R^3$ the body force per unit volume. The infinitesimal strain tensor, is defined as the symmetric part of the displacement gradients:
\begin{equation}\label{eq-strain}
\epsilon_{ij}({\bf u}) = \frac{1}{2}\left(u_{i,j} + u_{j,i}\right),
\end{equation}
where $u_{i,j} = \partial u_i/\partial x_j$.
The strain tensor is connected to the Cauchy stress tensor $\sigma_{ij}$ via the generalized Hooke's law:
\begin{equation}\label{eq-hooke}
\sigma_{ij}({\bf u}) = C_{ijkl}({\bf x})\epsilon_{kl}({\bf u}).
\end{equation}
where the material tensor $C_{ijkl}({\bf x})$ is a symmetric, positive definite fourth order tensor. The studied material will be further described in Section \ref{sec:compositeDescription}.

Now, let $\Gamma_{D}$ and $\Gamma_{N}$ be disjoint open subsets of $\Gamma$ such that $\overline{\Gamma_{D}}\cup\overline{\Gamma_{N}}=\Gamma$ and consider the function space \begin{equation}\label{eq-fctspace}
V := \{{\bf v} \in H^1(\Omega;\mathbb{R}^{3}) : {\bf v} = {\bf 0} \quad \text{on} \;\;\Gamma_D \}.
\end{equation}
Given functions ${\bf h}: \Gamma_D \rightarrow {\mathbb R}^{3}$ and ${\bf g}: \Gamma_N \rightarrow {\mathbb R}^{3}$, prescribing the Dirichlet and Neumann boundary data, the weak formulation of the problem to be considered consists in seeking the unknown displacement field ${\bf u} \in H^1(\Omega;\mathbb{R}^{3})$ with ${\bf u}={\bf h}$ on $\Gamma_{D}$ such that
\begin{equation}\label{eq-fe}
a({\bf u}, {\bf v}) = b({\bf v}) \quad \forall \mathbf{v} \in V,
\end{equation}
where the bilinear form $a(\cdot,\cdot)$ and the functional $b(\cdot)$ are defined by
\begin{equation}
\label{weak_form}
    a({\bf u}, {\bf v}) = \int_\Omega \sum_{i,j}\sigma_{ij}({\bf u})\epsilon_{ij}({\bf v}) \,dx \quad \text{and} \quad  b({\bf v})=\int_{\Gamma_{N}} \mathbf{g}\cdot{\bf v}\,ds + \int_\Omega \mathbf{f}\cdot\mathbf{v}\,dx.
\end{equation}

The variational problem \eqref{eq-fe} is discretized with conforming FEs on a mesh $\mathcal T_h$ on $\Omega$ by introducing the FE space $V_h \subset V$ as the tensor product $V_h \coloneqq V_h^1 \times V_h^2 \times V_h^3$
of the spaces $V_h^j := \text{span}\{\phi_{j}^{(i)}\}_{i=1}^N$, spanned by the usual Lagrange bases on $\mathcal T_h$. We find a function ${\bf u}^{p}\in H^1(\Omega;\mathbb{R}^{3})$ such that ${\bf u}^{p} = {\bf h}$ on $\Gamma_{D}$ and then seek an approximation ${\bf u}_h = {\bf u}^{p} + \widetilde{\bf u}_{h}$, where $\widetilde{\bf u}_{h}\in V_{h}$, such that
\begin{equation}\label{eqn:bilinear}
a(\widetilde{\bf u}_h,{\bf v}_h) = \widetilde{b}({\bf v}_h):=b({\bf v}_h)-a({\bf u}^{p}, {\bf v}_h) \quad \forall{\bf v}_h \in V_h.
\end{equation}
We block together displacements from all three space dimensions, so that ${\bf u}^{(i)}_h \in \mathbb B := \mathbb R^3$ denotes the vector of displacement coefficients containing all space components associated with the $i^{th}$ basis function. The displacement vector at a point $\bf x$ is then given by ${\bf u}_{h,j} ({\bf x}) = \sum_{i=1}^N {u}_{h,j}^{(i)} \;\boldsymbol \phi_{j}^{(i)}({\bf x})$, $j\in\{1,2,3\}$. The system \eqref{eqn:bilinear} is equivalent to a symmetric positive-definite (spd) system of algebraic equations:
\begin{equation}\label{eqn:fe_matrix_system}
{\bf A}{\bf \tilde{u}} = {\bf \tilde{b}} \quad \mbox{where} \quad {\bf A} \in \mathbb R^{N \times N} \quad \mbox{and} \quad {\bf \tilde{b}} \in \mathbb R^N.
\end{equation}
The blocks in the  global stiffness matrix and in the load vector, for any $i, j = 1, \ldots, N$, are given by ${\bf A}_{ij} = a(\boldsymbol \phi^{(i)}, \boldsymbol \phi^{(j)})$ and ${\bf \tilde{b}}_i = \widetilde{b}(\boldsymbol \phi_i)$. The vector ${\bf \tilde{u}} = [\widetilde{\bf u}_h^{(1)}, \ldots, \widetilde{\bf u}_h^{(N)}]^T \in \mathbb B^N$ is the block vector of unknown FE coefficients.

System \eqref{eqn:fe_matrix_system} can be assembled elementwise using Gaussian integration:
\begin{equation}\label{eq:element_stiffness}
a(v,w) = \sum_{e \in \mathcal T_h} a_e(v|_e, w|_e) \quad \forall v,w \in V. 
\end{equation}
The elementwise bilinear form $a_e$ is trivial to obtain here by restricting the integrals in \cref{eq-fe} to $e$. Later, restrictions of $a(\cdot,\cdot)$ to mesh-resolved subdomains will be crucial in defining coarse space approximations. For any mesh-resolved subdomain $D \subset \Omega$, the restriction of $a(\cdot,\cdot)$ to $D$ is denoted by

\begin{equation}\label{eqn:restricted_bilinear_form}
a_D(v,w) \coloneqq \sum_{e \in D} a_e(v|_e, w|_e) \quad \forall v,w \in V.
\end{equation}

\section{Multiscale-spectral generalized finite element method}

In this section, the methods employed in this paper are defined alongside some theoretical results. In addition to generalized finite element methods (GFEM), we also describe the \gls{GenEO} space, which was originally designed as a coarse space for two-level additive Schwarz methods in \cite{spillane2014abstract}. Here, it is for the first time applied as a multiscale method for stand-alone coarse approximation to a realistic three-dimensional multiscale problem in composites, in the form of the \gls{MS-GFEM} method \cite{babuska2011optimal,ma2022novel,ma2021error}.

\subsection{Domain decomposition}\label{sec:dd}

\begin{figure*}
    \centering
    \includegraphics[width=\linewidth]{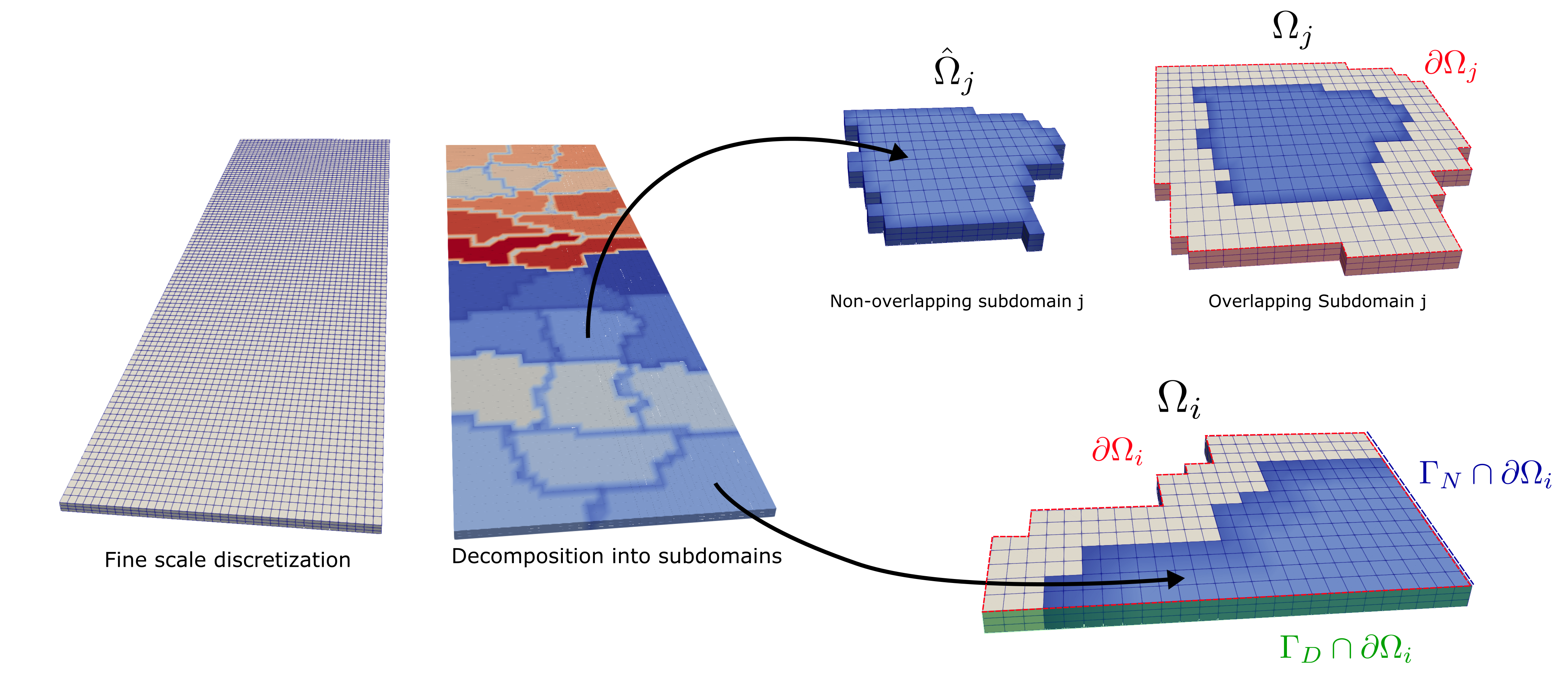}
    \caption{Domain decomposition into $24$ subdomains, illustrating the various components.}
    \label{fig:subdomain_terminology}
\end{figure*}
Both the two-level additive Schwarz method and GFEM are based on a decomposition of the domain $\Omega$ into non-overlapping subdomains $\{\hat{\Omega}_j\}_{j=1}^N$ that are resolved by the mesh $\mathcal{T}_h$; see \Cref{fig:subdomain_terminology} for an example. Each non-overlapping subdomain $\hat{\Omega}_j$ is extended by adding layers of neighboring elements to create an overlapping partition $\{\Omega_j \}_{j=1}^{N}$ of $\Omega$.

Next, the local \gls{FE} spaces 
\[
V_{h}(\Omega_j) \coloneqq \{v |_{\Omega_j} : v \in V_h\} \quad \text{and} \quad
V_{h,0}(\Omega_j) \coloneqq \{v |_{\Omega_j} : v \in V_h, \text{supp}(v) \subset \Omega_j\},
\]
are defined, where the former restricts $V_h$ to the subdomain $\Omega_j$ while the latter restricts this space further to functions whose support is contained entirely in $\Omega_j$.

A key ingredient of \gls{GenEO}-type coarse spaces and of GFEM is a \gls{PoU} subordinate to the overlapping decomposition $\{\Omega_j \}_{j=1}^{N}$. A particular partition of unity specific to the FE setting was constructed in \cite{spillane2014abstract}. For each $1\leq j\leq N$, let
\begin{equation}
    \text{dof}(\Omega_j) = \{k:1\leq k\leq n \,, \;\text{supp}({\boldsymbol \phi}^{(k)})\subset \overline{\Omega_j}\}
\end{equation}
denote the set of internal degrees of freedom in $\Omega_j$, and define for each degree of freedom $k \in \text{dof}(\Omega_j)$ a weight $\mu_{j, k}\in [0, 1]$ such that
\[ 
\sum_{\{j: 1 \leq j \leq N, k \in \text{dof}(\Omega_j)\}}  \mu_{j, k} = 1. 
\]
With these weights we can define a family of local partition of unity operators $\Xi_j: V_{h}(\Omega_j)\rightarrow  V_{h,0}(\Omega_j)$, $1\leq j\leq N$, such that
\begin{equation}\label{eqn:PoU}
\Xi_j ( v) \coloneqq \sum_{k \in \text{dof} ( \Omega_j)}  \mu_{j,k} v_k {{\boldsymbol \phi}^{(k)}}|_{\Omega_j},
\text{ for any } v=\sum_{k=1}^{n}v_{k}{\boldsymbol \phi}^{(k)} \in V_h ( \Omega_j).
\end{equation}
It follows immediately from the definition that the operators $\{ \Xi_{j}\}_{j=1}^{N}$ satisfy
\begin{equation}
  \sum_{j=1}^{N}R_{j}^{\top}\Xi_{j}(v|_{\Omega_j}) = v\qquad \text{for any} \;\;v\in V_{h}.
\end{equation}
Here $R_j^\top : V_{h,0}(\Omega_j) \rightarrow V_h $ denotes the prolongation operator defined as the extension of a function in $V_{h,0}(\Omega_j)$ by zero.
In \cite{spillane2014abstract}, it was suggested that the weights can be set as
\[ 
\mu_{j,k} \coloneqq \frac1{\#  \{ i : 1 \leq i \leq N\,,\; k \in \text{dof}(\Omega_i) \}}\;, 
\]
that is, one over the number of subdomains that contain $k$ as an internal \gls{DoF}. Other partitions of unity can also be used. In the numerical experiments below, we will use a different  (smoother) partition of unity.

\subsection{\gls{GenEO} coarse space}

The \gls{GenEO} space -- designed in the context of additive Schwarz preconditioning methods, as a robust coarse space correction for multiscale variational problems  -- is based on the following \gls{GEVP} on each subdomain $\Omega_{j}$:
Find $\lambda^{j} \in \mathbb{R}$, $\varphi^{j}_h \in V_{h}(\Omega_j)$ such that
\begin{equation}\label{eqn:bilinearGeneralizeEigenProblem}
a_{\Omega_j}(\varphi^{j}_h, v_h) = \lambda^{j}\,a_{\Omega_j}(\Xi_j ({\bf \varphi}^{j}_h), \Xi_j (v_h)), \quad \mbox{for all} \quad {v}_h \in V_{h}(\Omega_j).
\end{equation}
In the original publication \cite{spillane2014abstract}, the bilinear form on the right hand side of \eqref{eqn:bilinearGeneralizeEigenProblem} was restricted to the overlap $\Omega_j \backslash \hat{\Omega}_j$, but as shown in subsequent publications, such as \cite{bastian2021multilevel}, the GenEO space defined by the GEVP in \eqref{eqn:bilinearGeneralizeEigenProblem} has very similar coarse space correction properties. 

Only the lowest-energy eigenfunctions in \eqref{eqn:bilinearGeneralizeEigenProblem}, i.e., the ones corresponding to the smallest eigenvalues, are used to define the \gls{GenEO} coarse space.  
Denote by $\lambda^{j,k}$ and $\varphi_h^{j,k}$ the $k$-th smallest eigenvalue and the corresponding eigenfunction on subdomain $\Omega_j$. Then, the \gls{GenEO} coarse space is defined by 
\begin{equation}\label{eq:coarse_space}
 V_H \coloneqq \text{span} \Big\{ R_j^\top \Xi_j ( \varphi_h^{j,k} ) : k = 1, \ldots, m_j, \quad j = 1, \ldots, N \Big\},
\end{equation}
where the partition of unity operators are used to "stitch" the local approximation spaces on the subdomains $\Omega_j$ together and to guarantee that
$V_H \subset V_h$.
This definition of $V_H$ still leaves open the number of eigenfunctions $m_j$ to be included. 

In the context of two-level additive Schwarz methods, where the GenEO coarse space is combined (additively) with local solves on the overlapping subdomains $\Omega_j$ to obtain a preconditioning matrix ${\bf M}$ for ${\bf A}$, it is then possible to bound the condition number $\kappa$ of the preconditioned system independently of the mesh size $h$, the subdomain size $H$ or the heterogeneity in the coefficient. In particular, \cite[Corollary 3.23]{spillane2014abstract} states
\begin{equation}
\label{eq:condition_bound}
    \kappa({\bf M}^{-1} {\bf A}) \leq C(k_{0}) \max_{1 \leq j \leq N} \left( 1 + \frac{1}{\lambda^{j,m_j+1}} \right),
\end{equation}
where $C(k_0)$ is a typically small constant depending only on $k_{0}$, the maximum number of subdomains overlapping at any point. Thus, the condition number can be controlled by choosing the number $m_j$ of eigenfunctions per subdomain such that $1/\lambda^{j,m_j+1}$ is bounded uniformly across all subdomains. The eigenvalues $\lambda^{j,m_j+1}$ converge to $1$ as $m_j$ increases, but no theoretical results on the rate of convergence exist in general.

\subsection{Generalized FE methods with \gls{GenEO}-type local approximation}\label{sec:geneo_as_gmsfem}

The GenEO space $V_H$ in \eqref{eq:coarse_space} is in essence a global approximation space of generalized FEM type \cite{melenk}, where a fairly arbitrary family of local approximation spaces can be "stitched" together via a partition of unity to build the global space. As such, for $H$ (the subdomain size) sufficiently small or for $m_j$ sufficiently large, it is possible to solve the \gls{FE} problem directly to a required accuracy in $V_H$, in the spirit of the GFEM. However, as we will see below in the numerical experiments, the rate of convergence with respect to $m_j$ is rather poor when the local bases are computed as in \eqref{eqn:bilinearGeneralizeEigenProblem}. A significantly more efficient GFEM can be designed by slightly modifying the GEVP \eqref{eqn:bilinearGeneralizeEigenProblem} as shown in the following.

In this subsection, a particular family of \gls{GenEO}-type local approximation spaces is constructed and used within the framework of the GFEM as a stand-alone coarse approximation. Compared with the original version, there are two key ingredients in this \gls{GenEO}-type coarse space that provide a better accuracy for coarse approximation. The first is oversampling. Similarly to the construction of the overlapping subdomains, we extend each overlapping subdomain $\Omega_j$ further by adding more layers of fine-mesh elements to create an oversampling subdomain $\Omega_{j}^{\ast}$. The local eigenproblems used for constructing the coarse space will be defined on $\Omega_{j}^{\ast}$ instead of $\Omega_{j}$. The second ingredient is A-harmonicity. To make this notion precise, we first introduce the following local FE spaces defined on the oversampling domains:
 \begin{equation}\label{eqn:HD}
V_{h,{D}}(\Omega^{\ast}_j) = \{v \in V_h(\Omega^{\ast}_j) : v = 0 \text{ on } \partial\Omega^{\ast}_j \cap \Gamma_D\},
\end{equation}
\begin{equation}\label{eqn:HDI}
V_{h,{DI}}(\Omega^{\ast}_j) = \{v \in V_h(\Omega^{\ast}_j) : v = 0 \text{ on } \partial\Omega^{\ast}_j \cap (\Gamma_D \cup\Omega)\}.
\end{equation}
The space $V_{h,{D}}(\Omega^{\ast}_j)$ consists of FE functions restricted to $\Omega_j^{\ast}$ that vanish on the external Dirichlet boundary of $\Omega_j^{\ast}$, whereas $V_{h,{DI}}(\Omega_j)$ consists of FE functions that vanish on both the external Dirichlet boundary and the interior boundary of $\Omega_j^{\ast}$. The A-harmonic local \gls{FE} space on $\Omega^{\ast}_j$ is then defined as
\begin{equation}\label{eqn:AharmonicCondition}
V_{A}(\Omega^{\ast}_j) = \{u \in V_{h,{D}}(\Omega^{\ast}_j) : a_{\Omega^{\ast}_j}(u, v) = 0 \quad \forall v \in V_{h,{DI}}(\Omega^{\ast}_j)\}.
\end{equation}
Functions in $V_{A}(\Omega^{\ast}_j)$ are referred to as A-harmonic \gls{FE} functions. As we will see below, the local eigenvectors used for building the coarse space are A-harmonic \gls{FE} functions instead of general \gls{FE} functions. 

With the above notations, we now define a local eigenproblem similar to \eqref{eqn:bilinearGeneralizeEigenProblem} on each oversampling subdomain: Find $\lambda^{j} \in \mathbb{R}$, $\varphi^{j}_h \in V_{A}(\Omega^{\ast}_j)$ such that
\begin{equation}\label{eqn:AHarmonicGeneralizeEigenProblem}
a_{\Omega^{\ast}_j}(\varphi^{j}_h, v_h) = \lambda^{j}\,a_{\Omega^{\ast}_j}\big(\Xi_j ({\bf \varphi}^{j}_h|_{\Omega_j}), \Xi_j (v_h|_{\Omega_j})\big), \quad \mbox{for all} \quad {v}_h \in V_{A}(\Omega^{\ast}_j).
\end{equation}
Note that since $\Xi_j ({\bf \varphi}^{j}_h|_{\Omega_j})$ and $\Xi_j (v_h|_{\Omega_j})$ can be identified with \gls{FE} functions in $V_{h}(\Omega_j^{\ast})$, the right-hand side of the above \gls{GEVP} is well-defined.

Let $(\lambda^{j,k},\,\varphi^{j,k}_h)$ denote the $k$-th eigenpair of the GEVP \eqref{eqn:AHarmonicGeneralizeEigenProblem} with eigenvalues enumerated in increasing order. The desired \gls{GenEO}-type GFEM coarse space is defined almost identically to the standard version \eqref{eq:coarse_space}:
\begin{equation}\label{eq:AHarmonic_coarse_space}
 V_H \coloneqq \text{span} \Big\{ R_j^\top \Xi_j ( \varphi_h^{j,k}|_{\Omega_j} ) : k = 1, \ldots, m_j, \quad j = 1, \ldots, N \Big\}.
\end{equation}

The last ingredient of the \gls{MS-GFEM} method is a global particular function built from local particular functions. On each oversampling subdomain $\Omega_{j}^{\ast}$, we first define a local particular function ${\bf u}^{p}_{h,j} = \psi^{r}_{h,j}+\psi^{d}_{h,j}$, where $\psi^{r}_{h,j}\in V_{h,DI}(\Omega_{j}^{\ast})$ satisfies
\begin{equation}\label{eq:up_r}
a_{\Omega^{\ast}_j}(\psi^{r}_{h,j}, v_h) = b_{\Omega_j^{\ast}}(v_{h})\quad \forall v_{h}\in V_{h,DI}(\Omega_{j}^{\ast})
\end{equation}
with $b_{\Omega_j^{\ast}}(\cdot)$ being the restriction of $b(\cdot)$ to $\Omega^{\ast}_{j}$, and $\psi^{d}_{h,j}\in V_{h}(\Omega_j^{\ast})$ satisfies $\psi^{d}_{h,j}={\bf h}$ on $\Gamma_{D}\cap \partial \Omega_{j}^{\ast}$ and
\begin{equation}\label{eq:up_d}
a_{\Omega^{\ast}_j}(\psi^{d}_{h,j}, v_h) = 0\quad \forall v_{h}\in V_{h,D}(\Omega_{j}^{\ast}).
\end{equation}
Note that $\psi^{d}_{h,j}$ vanishes on all interior subdomains where $\Gamma_{D}\cap \partial \Omega_{j}^{\ast} =\emptyset$ or whenever ${\bf h}={\bf 0}$ on $\Gamma_{D}\cap \partial \Omega_{j}^{\ast}$. On subdomains intersecting $\Gamma_{D}$ it would in fact be possible to combine problems \eqref{eq:up_r} and \eqref{eq:up_d} into one local problem, but this leads to a slightly larger constant $C$ in Theorem~\ref{thm:error_bound} below. Therefore, we work with local particular functions defined via \eqref{eq:up_r} and \eqref{eq:up_d} in this paper.
The global particular function is then defined by ``stitching'' together the local functions using the partition of unity:
\begin{equation}\label{eq:global_particular_function}
    {\bf u}_{h}^{p}=\sum_{j=1}^{N} R_j^\top \Xi_j ( {\bf u}^{p}_{h,j}|_{\Omega_j} ).
\end{equation}

Having defined the coarse space $V_{H}$ and the global particular function ${\bf u}_{h}^{p}$, we are now ready to give the MS-GFEM method for solving the fine-scale \gls{FE} problem \eqref{eqn:bilinear}: Find ${\bf u}_{h}^{G} = {\bf u}_{h}^{p} + {\bf u}_{H}$, where ${\bf u}_{H}\in V_{H}$, such that
\begin{equation}\label{eq:GFEM_coarse_approximation}
 a({\bf u}_{h}^{G}, {\bf v}) = b({\bf v})\quad \forall {\bf v}\in V_{H}.   
\end{equation}

To assess the quality of the MS-GFEM approximation, we estimate the energy norm $\Vert {\bf v}\Vert_{a}:=\sqrt{a({\bf v},{\bf v})}$ of the error ${\bf u}_h - {\bf u}_{h}^{G}$. Following the lines of the proofs of Theorems 2.1 and 3.4 in \cite{ma2022novel}, we can derive the following bound:
\begin{theorem}\label{thm:error_bound}
\begin{equation}\label{eq:error_bound}
    \| {\bf u}_h - {\bf u}_{h}^{G} \|_a \leq 
    C \left(\lambda^{\not\in}_{\min} \right)^{-1/2} \|{\bf u}_h\|_a \quad \text{where} \quad \lambda^{\not\in}_{\min} :=
    \min_{1 \leq j \leq N} 
    \lambda^{j,m_j+1} , 
\end{equation}
$C$ is known explicitly and bounded by the maximum number of oversampling domains that overlap at 
any given point in $\Omega$,
and $\lambda^{j,m_j+1}$ is the smallest eigenvalue corresponding to any eigenvector not included $(\not\in)$ in the local basis on $\Omega_j$.
\end{theorem}

Thus, the efficiency of the MS-GFEM method is controlled by the speed at which the eigenvalues in \eqref{eqn:AHarmonicGeneralizeEigenProblem} grow. 
To estimate this growth rate, let $H_j$ and $H_{j}^{\ast}$ denote the diameter of $\Omega_j$ and $\Omega_{j}^{\ast}$, respectively. 
The following (informal) theorem in $d$ dimensions, similar to Theorem 4.6 in \cite{ma2021error}, gives an exponential bound on the eigenvalues, which can be proved following the lines of the proofs of Theorem 4.6 in \cite{ma2021error} and Theorem 7.3 in \cite{babuvska2014machine} (details will be given in a forthcoming paper).
\begin{theorem}\label{thm:exponential_decay}
Let $h$ be sufficiently small. Then there exist 
$k_{j}, b_j, C_j >0$ independent of $h$, such that for $k>k_{j}$
\begin{equation}\label{eq:decay_rate}
\left(\lambda^{j,k}\right)^{-1} \leq C_j e^{-{b_j}k^{{1}/{(d+1)}}}.
\end{equation}
The constants $C_j$, $k_j$ and $b_j$ can again be derived explicitly. The value of $b_j$ and thus the convergence rate grows with the amount of oversampling, i.e., with decreasing $H_{j}/H_{j}^{\ast}$.
\end{theorem}

Combining the exponential bound \eqref{eq:decay_rate} on the local eigenvalues and the global error estimate \eqref{eq:error_bound} provides a rigorous, exponential error bound for the MS-GFEM method. It is important to note that the exponential growth rate of the local eigenvalues critically relies on the two aforementioned ingredients of the new coarse space, i.e., oversampling and A-harmonicity. The global error estimate \eqref{eq:error_bound} also holds for the standard \gls{GenEO} coarse space when used as a coarse approximation. However, without the two key ingredients, the eigenvalues of the local \gls{GEVP}~\eqref{eqn:bilinearGeneralizeEigenProblem} do not grow exponentially fast, 
making the standard \gls{GenEO} coarse space significantly less efficient; see Subsection~\ref{sec:comparison_two_coarse spaces}.

Apart from the exponential decay rate of the error with respect to the number of local basis functions, Theorem~\ref{thm:exponential_decay} also provides an explicit decay rate of the error with respect to the oversampling size, which offers a second handle to control the error of the method, in addition to a change in the size of the local approximation spaces. This turns out to be of great importance in reducing the size of the global coarse problem; see Subsection~\ref{sec:scalability}.

We end this subsection by discussing ways to solve the local GEVP \eqref{eqn:AHarmonicGeneralizeEigenProblem}. Due to the presence of the A-harmonic condition, a straightforward yet time-consuming way of solving \eqref{eqn:AHarmonicGeneralizeEigenProblem} is to first construct the basis functions of the A-harmonic FE space by solving many local boundary value problems \cite{babuska2011optimal,babuska2020multiscale,chen2020randomized}. Instead, we use a different and more efficient method proposed in \cite{ma2021error}, where the A-harmonic condition is directly incorporated into the local \gls{GEVP}. To this end, a Lagrange multiplier is introduced and the local \gls{GEVP} \eqref{eqn:AHarmonicGeneralizeEigenProblem} is rewritten in an equivalent mixed formulation: Find $\lambda \in \mathbb{R}$, $\varphi_h \in V_{h,D}(\Omega^{\ast}_j)$ and $p_h \in V_{h,{DI}}(\Omega^{\ast}_j)$ such that
\begin{equation}\label{eqn:bilinearAharmonic}
\begin{aligned}
a_{\Omega^{\ast}_j}(\varphi_h, v_h) + a_{\Omega^{\ast}_j}( v_h, p_h) &= \lambda \,a_{{\Omega}^{\ast}_j}(\Xi_j ({\bf \varphi}_h|_{\Omega_j}), \Xi_j(v_h|_{\Omega_j})),& \quad \forall & {v}_h \in V_{h,D}(\Omega^{\ast}_j),\\
a_{\Omega^{\ast}_j}( \varphi_h, \xi_h) &= 0, &\quad \forall &{\xi}_h \in V_{h,{DI}}(\Omega^{\ast}_j).
\end{aligned}
\end{equation}
The implementation details of how the augmented system \eqref{eqn:bilinearAharmonic} is solved are presented in Subsection~\ref{sec:A-harmo_implementation}.

\section{Implementation}


The \gls{DUNE} \cite{DUNE} package is an open-source, modular toolbox for the numerical solution of PDE problems. It leverages advanced C++ programming techniques in order to provide modularity from the ground up while producing highly efficient applications. As such, it allows the reuse of many existing components when implementing our new mathematical methods for \gls{HPC} applications.
The new methods were integrated in the \texttt{dune-composites} module \cite{reinarz2018dune}, which facilitates setting up elasticity models and provides access to efficient solvers that scale to thousands of cores on modern \gls{HPC} systems, despite the typically bad conditioning of composites problems. This was achieved through an \gls{HPC}-scale GenEO implementation \cite{geneo-hpc} developed as part of \texttt{dune-composites} and later moved into the lower-level discretization module \texttt{dune-pdelab} \cite{PDELAB} within \gls{DUNE}. The \texttt{dune-pdelab} module and several lower-level \gls{DUNE} modules  
are used within \texttt{dune-composites} to obtain the finite element discretizations on the fine level.

Within the \gls{DUNE} framework a number of grid implementations are provided for various purposes. Initially, only \texttt{YASPGrid}, a structured grid, provided native support for overlapping domain decomposition methods, 
where each process holds a copy of elements in its overlap region and may exchange data attached to associated \gls{DoF}s with neighboring processes. Since the methods considered here are
constructed using overlapping subdomains, we need this kind of communication mechanism. In previous work within \texttt{dune-composites} \cite{reinarz2018dune}, the restriction to cuboid domains induced by \texttt{YASPGrid} was overcome by applying a geometric transformation to the grid. However, for many relevant engineering applications, a smooth transformation from a cuboid geometry to the actual geometry is not available. Thus, \gls{DoF}-based communication was extended to support unstructured grids as well. In order to handle more general model geometries, the grid creation was first shifted to a pre-processing step using \emph{Gmsh}~\cite{geuzaine2009gmsh}. The grid import is facilitated by the IO functionality of \texttt{dune-grid}. 

\subsection{Algebraic overlap construction for GenEO}

In order to provide overlap and communication across overlaps on unstructured grids, an algebraic approach was introduced where finite element matrices $A_j$ defined on overlapping subdomains $\Omega_j$ are constructed from corresponding matrices $\hat{A}_j$ assembled on non-overlapping subdomains $\hat{\Omega}_j$. 

\Cref{fig:subdomain_terminology} illustrates our notation for non-overlapping and overlapping subdomains.
We begin with a matrix $\hat{A}_j$ on a non-overlapping subdomain $\hat{\Omega}_j$, and assume it has neighboring subdomains $\hat{\Omega}_k$, $k\in E_j$, with corresponding matrices $\hat{A}_k$. We identify coinciding \gls{DoF}s in $\hat{A}_j$ and $\hat{A}_k$ through a global indexing, as provided by non-overlapping \gls{DUNE} grids. On subdomain $k$, we now identify all \gls{DoF}s directly connected to those in $\hat{\Omega}_j$, and provide process $j$ with unique indices and the connectivity graph of this newly identified layer of \gls{DoF}s. Since connected \gls{DoF}s in a \gls{FE} discretization are either associated with the same element or adjacent elements, we now are in a position where each process is aware of its neighbors' \gls{DoF}s within the first layer of elements along the respective boundary.

\begin{figure}[t]
    \centering
    \includegraphics[width=.8\linewidth]{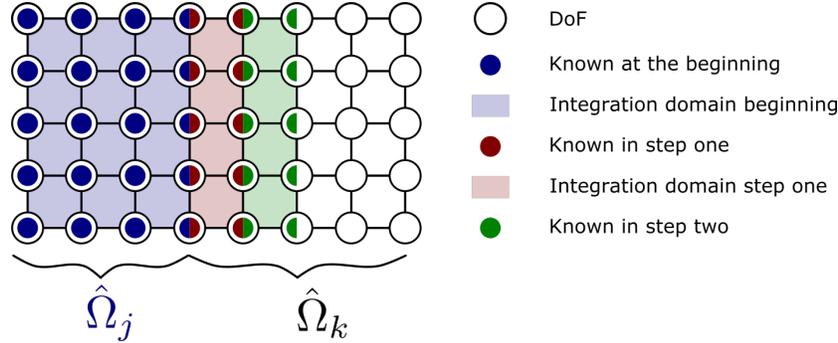}
    
    \caption{The degrees of freedom known to process $j$ after each step of recursive
    extension of the matrix connectivity graph.}
    \label{fig:extending_connectivity_graph}
\end{figure}

By a recursive application of this algorithm, as shown in \Cref{fig:extending_connectivity_graph}, we can now grow the algebraic overlap by an arbitrary number of layers of elements, without having to rely on the grid implementation to provide connectivity graphs.
In the following, we will define the number of layers of elements added to a subdomain as $o$ (for overlap) -- even though the overlap of a domain with its neighbour is in fact $2o$. 
Finally, using the extended connectivity graphs above we can construct a communication mechanism between neighbors to exchange overlap data of vectors defined on the extended domains. We now have the communication infrastructure to generate a \gls{GenEO} space in a parallel way in place .
While the method above allows us to extend connectivity graphs of non-overlapping matrices into neighboring subdomains, mimicking the connectivity graphs of the overlapping matrices $A_j$, the \gls{GenEO} method obviously requires the actual matrix entries. Using communication across the algebraic overlap, the required matrix entries can be exchanged between neighbors. This relies on the basic assumption (satisfied in general) that the bilinear form $a$ of the weak formulation can be decomposed additively into elementwise bilinear forms $a_{e}$, such that $a(u,v) = \sum_{e \in \mathcal{T}_h} a_{\epsilon}(u,v)$ with $\mathcal{T}_h$ the grid on $\Omega$. For linear elasticity this is trivially satisfied, see \eqref{eq:element_stiffness}.

\begin{figure}[t]
    \centering
    \includegraphics[width=.8\linewidth]{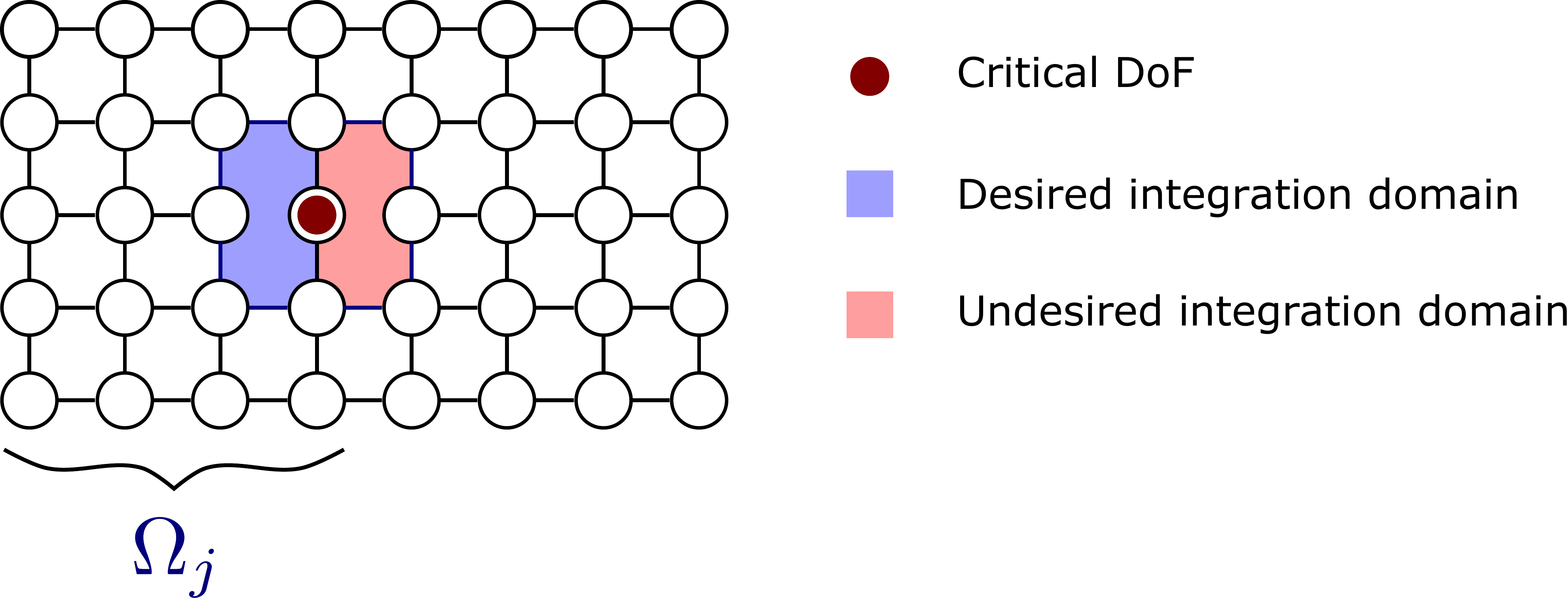}
    \caption{Integration domain of a \gls{DoF} on the boundary of the overlap of $\Omega_j$. Since it lies in the interior of a neighboring non-overlapping domain, its matrix entries are computed by the neighbor. The neighbor takes the entire integration domain around the \gls{DoF} into account, while process $j$ requires only integration up to the overlapping subdomain boundary $\partial \Omega_j$. The neighbors' matrix entries can therefore not be used to assemble $A_j$.}
    \label{fig:boundary_DoF_integration_domain}
\end{figure}

\begin{figure}[t]
    \centering
    \includegraphics[width=.6\linewidth]{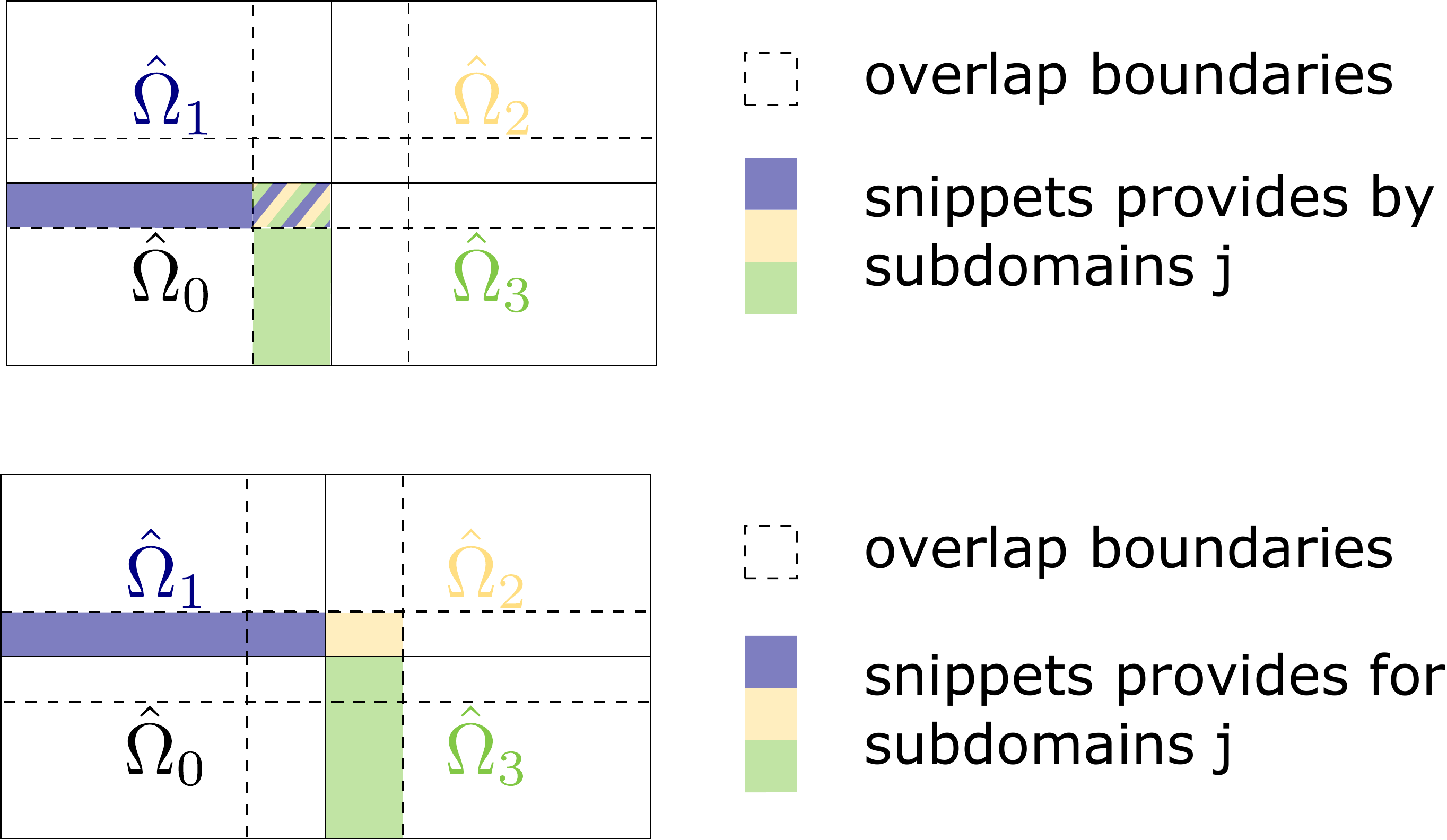}
    \caption{Subdomain snippets $\Omega_j \cap \hat{\Omega}_k$ used to assemble matrices on the algebraic overlaps, ensuring correct entries on overlap boundaries.}
    \label{fig:snippets}
\end{figure}

As a result, retrieving matrix entries from neighbors and adding them where \gls{DoF}s belong to multiple subdomains leads to a correct construction of an overlapping subdomain matrix from non-overlapping ones. There is one exception however: At the boundary $\partial \Omega_j$ of the overlapping domain, retrieving the entries of an interior \gls{DoF} will not yield the correct result since the integration domain of said interior \gls{DoF} extends beyond $\Omega_j$ as shown in \Cref{fig:boundary_DoF_integration_domain}.
We circumvent this issue by not sending individual matrix entries in $\hat{A}_k$ from process $k$ to process $j$, but instead by assembling a new matrix $\hat{A}^j_k$ on the domain $\Omega_j \cap \hat{\Omega}_k$ for that purpose as shown in \Cref{fig:snippets}. In practice, we communicate the partition of unity $\Xi_j$ generated from the already available matrix graph of $A_j$ to all neighbors, and use $\Omega_j \cap \hat{\Omega}_k = \text{supp}(\Xi_j) \cap \hat{\Omega}_k$ as a convenient proxy to the desired domain snippet. Since we have the decomposition 
\[ 
\Omega_j = \hat{\Omega}_j \cup \bigcup_{k\in E_j} \left( \Omega_j \cap \hat{\Omega}_k \right), 
\]
all snippets together form the desired overlapping subdomain. Exploiting additivity of the bilinear form $a$ and of the
\gls{FE} matrix entries, this assembly strategy delivers the correct overlapping matrix $A_j$ for $\Omega_j$.

\subsection{Partition of Unity and oversampling}

\begin{figure}
    \centering
    \includegraphics[width=.8\linewidth]{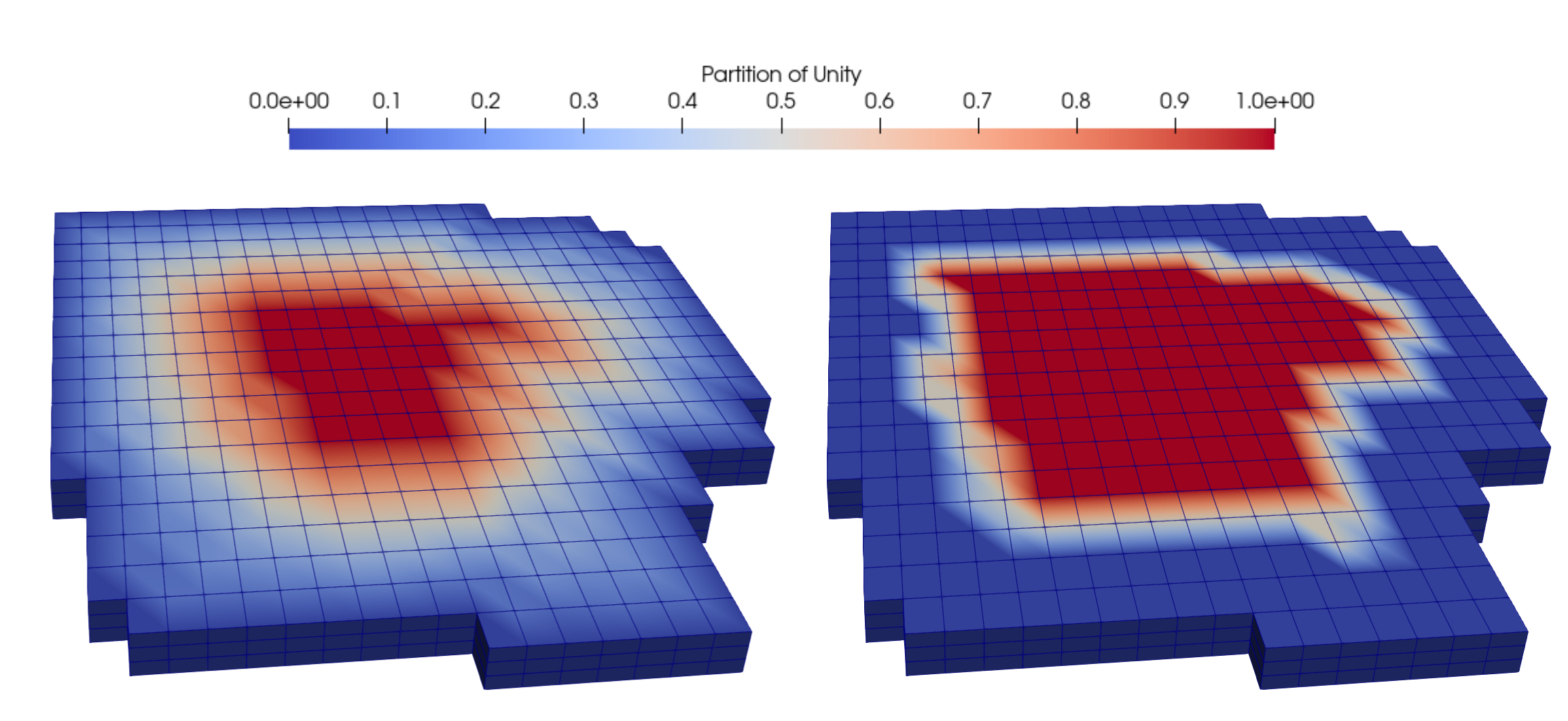}
    \caption{Partition of unity operators $\Xi_j$ for an overlap size of $o=3$ (left) and $o=1$ (right). On the right we also plot the oversampling domain $\Omega_{j}^{\ast}$ where $o^*=3$, as well as the partition of unity on $\Omega_j$ where $o=1$.}
    \label{fig:PoU-Oversampling}
\end{figure}

As explained in \cref{sec:dd}, various partition of unity operators can be used in the \gls{GenEO} theory, as long as they fulfill condition \eqref{eqn:PoU}. 
One suitable choice is a smooth transition from zero on $\partial\Omega_j$ to one on $\Omega_j \setminus \cup_{k \neq j} \Omega_k$, as illustrated in \Cref{fig:PoU-Oversampling} (left). It is the choice made for the \gls{GenEO} coarse space in \cite{reinarz2018dune} (where it is used as a preconditioner). 
To construct this partition of unity, a weight $w_p$ is attributed to each \gls{DoF} $p \in \text{dof}(\Omega_j)$.
Initially, $w_p$ is set to $0$ on $\partial\Omega_j$ and to $2o$ on the remaining vertices in $\Omega_j$, where $o$ is the number of layers added. For each \gls{DoF} $p$, using the extended connectivity graph of $\Omega_j$ in its vicinity and comparing with neighboring weights, the weight $w_p$ is then reduced incrementally in the subdomain interior such that
\[ 
w_p = \min\bigg(\min_{q_1 \leq q_i \leq q_n} (w_{p}, w_{q_i} + 1), 2o\bigg),
\]
where $\{q_{1} ,..., q_n \}$ are the neighboring \gls{DoF}s of $p$. It suffices to iterate this $2o-1$ times.
As an example, the partition of unity in \Cref{fig:PoU-Oversampling} (left) has been constructed with $o=3$. 
The communication mechanism between subdomains is used to associate to each \gls{DoF} $p \in \text{dof}(\Omega_j)$ their corresponding counterpart in the neighboring subdomain $\Omega_k$. 
The partition of unity is then simply defined as:
\[ 
\mu_{j,p} \in [0, 1]: \qquad \mu_{j, p} = \frac{w_{j, p}}{\sum_{\{1 \leq i \leq N\}} w_{i, p}}, \qquad p \in \text{dof}(\Omega_j)
\]
where $N$ is the number of subdomains sharing the \gls{DoF} $p$. 

The handling of oversampling subdomains is controlled via the choice of the partition of unity. 
The implementation of the partition of unity for the oversampled subdomains $\Omega_{j}^{\ast}$ is the same as the one described above with a different initialisation of $w_p$. In particular,
$w_p$ is initialized to $0$ not only on $\partial\Omega_j$, but also for a further $o^{\ast}-o$ layers of \gls{DoF}s before applying the iterative process above.
Thus, the partition of unity is a vector defined on the full oversampling subdomain $\Omega_{j}^{\ast}$ that takes the value zero on $\Omega_{j}^{\ast} \setminus \Omega_{j}$. In the following, we will choose a variable oversampling size $o^{\ast}$ and keep one layer of non-zero partition of unity overlap, i.e., $o=1$, as shown in \Cref{fig:PoU-Oversampling} (right).

\subsection{Enforcing A-harmonicity within the \gls{GEVP} in \gls{DUNE}}\label{sec:A-harmo_implementation}

As described in \cite{ma2021error}, the \cref{eqn:bilinearAharmonic} can be formulated as a matrix eigenvalue problem. To achieve that, the \gls{DoF} associate with $\Omega_{j}^{\ast}$ are partition into three sets:
\[
\mathcal{B}_1 = \text{dof}\Big((\Omega_{j}^{\ast} \setminus \partial\Omega_{j}^{\ast})\cup (\partial\Omega_{j}^{\ast} \cap \Gamma_N)\Big), \ \ 
\mathcal{B}_2 = \text{dof}\Big(\partial\Omega_{j}^{\ast} \setminus (\Gamma_D \cup \Gamma_N)\Big), \ \  
\mathcal{B}_3 = \text{dof}\Big(\partial\Omega_{j}^{\ast} \cap\Gamma_D\Big).
\]
Those sets of \gls{DoF} are depicted in \Cref{fig:subdomain_terminology}. We also define $n_i$ as the sizes of the corresponding set $\mathcal{B}_i$. The \gls{GEVP} in matrix form is defined as follows: Find $\lambda_h \in \mathbb{R}$, $\tilde{\phi_j} = (\phi_{j,1}, \phi_{j,2}) \in \mathbb{R}^{n_1+n_2}$ and $p \in \mathbb{R}^{n_1}$:
\begin{equation}\label{eqn:algebraicAharmonic}
\begin{pmatrix} A_{j,11}&A_{j,12}&A_{j,11}\\
A_{j,21}&A_{j,22}&A_{j,21}\\
A_{j,11}&A_{j,12}&0 \end{pmatrix}
\begin{pmatrix} \phi_{j,1}\\\phi_{j,2}\\p \end{pmatrix}
= \lambda 
\begin{pmatrix} B_{j,11}&0&0\\
0&0&0\\
0&0&0 \end{pmatrix} 
\begin{pmatrix} \phi_{j,1}\\\phi_{j,2}\\p \end{pmatrix},
\end{equation}
where ${A}_{j,mn} = a_{\Omega_{j}^{\ast}}\big(\varphi_{k}, \varphi_{l}\big)_{k\in \mathcal{B}_m,l\in \mathcal{B}_n}$ and ${B}_{j,11} = a_{\Omega_{j}^{\ast}}\big(\Xi_j \varphi_{k}, \Xi_j \varphi_{l}\big)_{(k,l)\in \mathcal{B}_1}$. The blocks ${B}_{12}$, ${B}_{21}$, ${B}_{22}$ are zero since the partition of unity vanishes on $\mathcal{B}_2$. 
The oversampling, created via the choice of partition of unity, will also affect the the right hand side of  \cref{eqn:algebraicAharmonic}, such that not only the blocks ${B}_{12}$, ${B}_{21}$, ${B}_{22}$ but also all entries in $B_{j,11}$ corresponding to the oversampling region $\Omega_{j}^{\ast} \setminus \Omega_{j}$ will be zero (see \Cref{fig:PoU-Oversampling}).
The GEVP solution to construct the local basis is then $\phi_j = (\phi_{j,1}, \phi_{j,2}, \phi_{j,3}) \in \mathbb{R}^{n_1+n_2+n_3}$, where $\phi_{j,3}$ is a zero-vector corresponding to the \gls{DoF}s in $\mathcal{B}_3$, which combined with $(\phi_{j,1}, \phi_{j,2})$ provides an A-harmonic coefficient vector on all of $\Omega_{j}^{\ast}$.
Instead of generating each block $A_{j,nm}$ individually, it is extracted from $A_j$ using the sets of \gls{DoF}s $\mathcal{B}_1, \mathcal{B}_2, \mathcal{B}_3$. In practice, the matrix
\[
\begin{pmatrix} A_{j,11}&A_{j,12}&A_{j,11}\\
A_{j,21}&A_{j,22}&A_{j,21}\\
A_{j,11}&A_{j,12}&0 \end{pmatrix} = \begin{pmatrix} M_1&M_2^T\\
M_2&0\end{pmatrix}
\]
is built by block. The top left block $M_1$ is obtained by removing the rows and columns corresponding to $\mathcal{B}_3$ from $A_j$. In the elasticity case, Dirichlet boundary conditions are imposed in \gls{DUNE} by altering rows in the matrix $A_j$ -- one on  the diagonal, zero elsewhere -- so the set $\mathcal{B}_3$ can be easily detected. Then finally $M_1$ is obtained by
removing the rows and columns corresponding to $\mathcal{B}_2$. The block that is removed corresponds to $M_2$. \gls{DoF}s belonging to $\mathcal{B}_2$ are detected and saved during the overlap creation phase.

Once the basis is obtained, the next step in the construction of the coarse space $V_H$ is to multiply the eigenvectors by the partition of unity, see \cref{eq:coarse_space}. 
To improve the conditioning of the coarse space problem, a re-orthogonalization step via a Gram-Schmidt process is carried out following this multiplication.
Contrary to \cite{ma2021error}, we have implemented directly \cref{eqn:algebraicAharmonic}, since the simplification of \cref{eqn:algebraicAharmonic} proposed in \cite{ma2021error} for diffusion problems requires a special handling of subdomains that are only affected by Neumann boundaries $\Gamma_N$, which is more involved for linear elasticity in 3D.

\subsection{Boundary conditions at the coarse level}

The eigenvectors $ (\phi^k_j)_{k \in (1,m_j)} $ of the local \gls{GEVP}s have to be combined with local particular solutions if the body force $\bf f$ is nonzero or if $\Omega_{j}^{\ast}$ touches the global Dirichlet boundary and the boundary displacement $\bf h$ is nonzero.
A particular solution $\psi_j$ is obtained by solving the local \gls{PDE}
\[
{A_j}{\psi_j} = {b_j} \quad \mbox{where} \quad {A_j} \in \mathbb B^{N_j} \times \mathbb B^{N_j} \quad \mbox{and} \quad {b_j} \in \mathbb B^{N_j}
\]
The solution $\psi_j$ is then multiplied by the partition of unity and orthogonalized with respect to $ (\phi^k_j)_{k \in (1,m_j)} $ using again Gram-Schmidt. We denote the resulting vector by $\hat{\psi}_j$. This vector is then normalized and appended to the basis. 

At the coarse level, inhomogeneous Dirichlet boundary conditions are imposed by altering the coarse matrix $A^H$ and the coarse vector $b^H$. If $\ell_j$ denotes the index of the particular solution $\psi_j$ on subdomain $\Omega_{j}^{\ast}$ in the coarse system, then a Dirichlet boundary condition, imposed on $\Omega_{j}^{\ast}$ via the particular solution $\psi_j$, is enforced in the coarse system by setting:
\[
A^H_{\ell_j \ell_j} = 1, \quad A^H_{\ell_j k} = 0 \ \ \text{for all} \ \ \quad k \neq \ell_j, \ \  
\text{and} \quad b^H_{\ell_j} = \lVert \hat{\psi}_j \rVert.
\]

\subsection{Hardware and software}

For solving the local \gls{GEVP}, we use Arpack~\cite{lehoucq1998arpack} through the Arpack++ wrapper in symmetric shift-invert mode. As subdomain solver, we use UMFPack \cite{UMFPACK} in the eigenvalue solver.
The domain partition of $\Omega$ into non-overlapping subdomains $\hat{\Omega}_j$ is carried out by the graph partitioner \emph{ParMetis}~\cite{karypis1998multilevelk}. 

The numerical results have been carried out on the Hamilton HPC Service at Durham University. Its last version, called \emph{Hamilton8}, provides a total of 15,616 CPU cores, 36TB RAM and 1.9PB disk space. \emph{Hamilton8} is composed of 120 standard compute nodes, each with 128 CPU cores (2x AMD EPYC 7702), 256GB RAM and 400GB local SSD storage.

\section{Numerical Experiments: Performance tests for MS-GFEM on composite structures}

Throughout this paper, we assume a linear elastic behavior of the composite material. For aerospace applications, considering damage behavior is essential. Ultimately, our framework will be applied to model large displacement effects and the onset of failure. To achieve this, two requirements need to be validated. Firstly, it is essential to have an efficient linear elastic solver for a large-scale model to be used within any non-linear iteration. Secondly, the approximate solution needs to accurately represent relevant damage criteria, especially local extrema.

In this section, after a brief description of the composite models, we first analyse the output from all components of the MS-GFEM when applied to a composite beam. 
The method is then applied to a complex, aerospace part using a composite failure criterion to assess the accuracy.
Finally, the parallel scaling of the method is investigated, demonstrating the efficiency of the proposed MS-GFEM on large composite structures.

\subsection{Specification of the considered composite model problems}\label{sec:compositeDescription}

\begin{figure}[t!]
    \centering
    \includegraphics[width=0.7\linewidth]{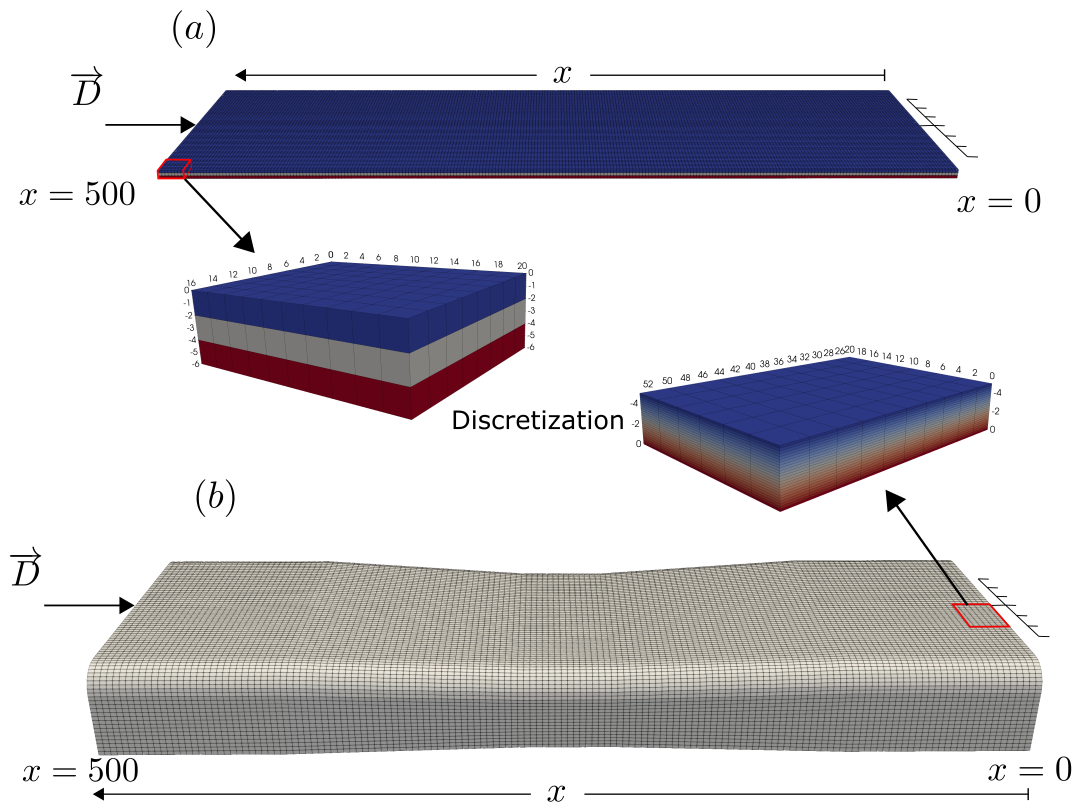}
    \caption{Models used in the numerical experiments: laminated beam (a) and laminated C-spar (b).\!\!}
    \label{fig:Models}
\end{figure}

In the following three subsections, first a laminated beam under compression is considered in order to evaluate the performance of the method proposed in this study. It is illustrated in \Cref{fig:Models} (a).
The laminated beam has a length of $500$mm, a width of $140$mm and a thickness of $6$mm. The laminate is made up of a stack of three layers (or plies) of the same thickness. Each layer represents a uni-directional composite made up of carbon fibres embedded into resin. Plies are modelled as homogeneous orthotropic elastic materials, characterised by nine parameters and a vector of orientations~${\bf\theta}$. The elastic properties of AS4/8552 (\cite{falco2018modelling}) have been chosen for this example. In the global coordinate system, the material tensor is orientated using standard tensor rotations, following the stacking sequence  [$0^{\circ}, 45^{\circ}, -45^{\circ}$]. For more details see, e.g., \cite{koay2009six}.
The study considers the elastic behavior of the laminated beam under compression. The uni-axial compression, illustrated in \Cref{fig:EValuesVectors}(b), is modeled as a displacement imposed Dirichlet boundary condition.

Then for the remainder of the paper, a realistic aerospace part is used to demonstrate the high quality of the achieved coarse approximation. 
The aerospace part in question is a $500$mm long C-shaped wing spar section (C-spar) with a joggle region in its center, creating a geometric feature in the structure, see \Cref{fig:Models} (b).
The material is a laminated composite, composed of 24 uni-directional layers (carbon fibers and resin) of $0.2mm$ each, which are orientated as follows:
\[
\big[(45^\circ, -45^\circ)_3, (0^\circ, 90^\circ)_3\big]_s
\]
The behavior of the C-spar under compression will be investigated using the same boundary conditions as the ones used in the beam example above.

In both examples, as visualized in \Cref{fig:Models}, for the local approximation we use FE grids with piecewise linear elements, reduced integration to avoid shear-locking and one element through thickness per layer.

\subsection{Local \gls{GEVP} outputs}

\begin{figure*}[t!]
    \centering
    \includegraphics[width=\linewidth]{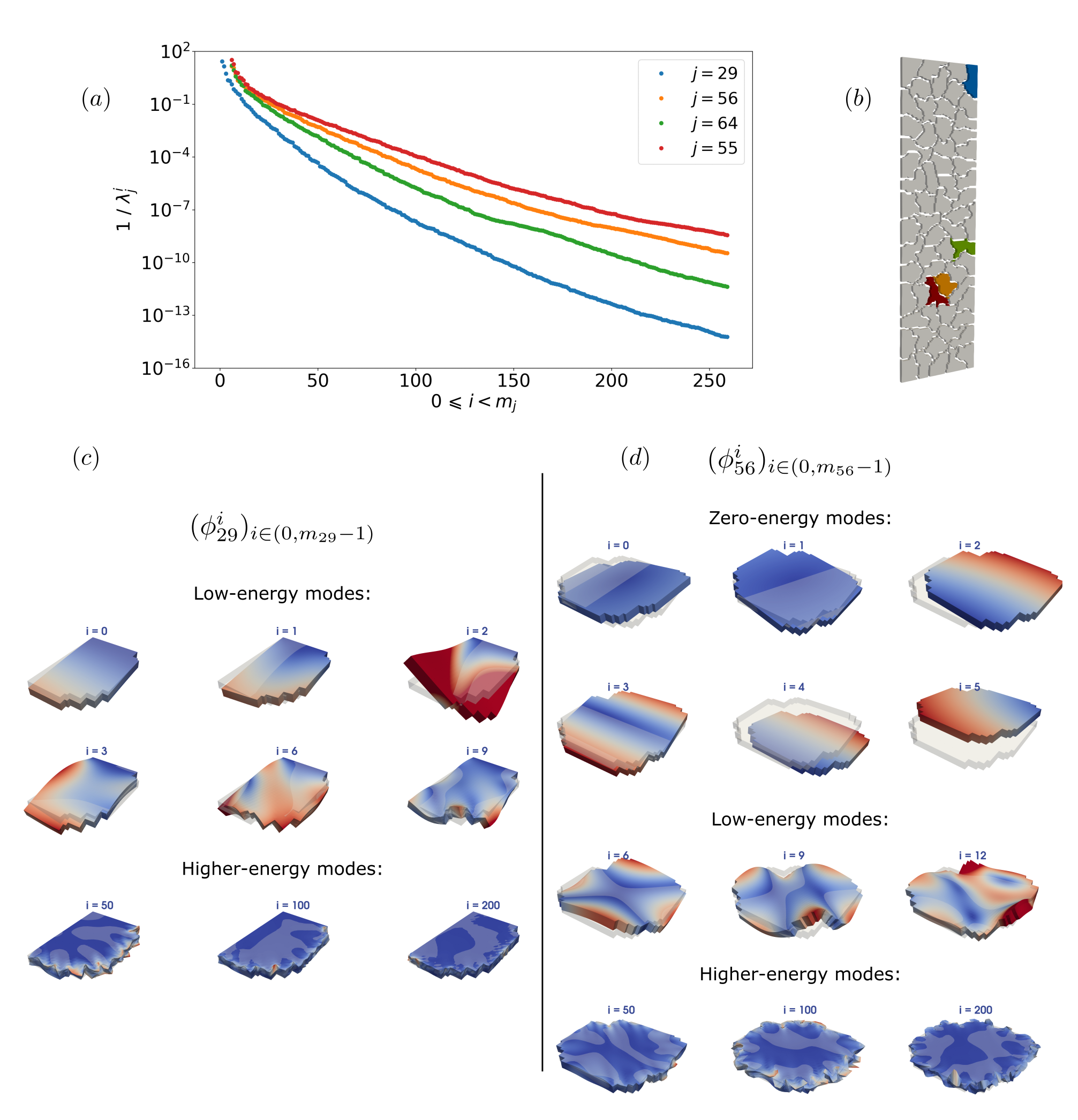}
    \caption{Visualisation of eigenpairs of the A-harmonic local \gls{GEVP} with oversampling $o^{\ast}=8$, i.e., eight layers of elements around subdomains. The reciprocal of the $i^{th}$ eigenvalue $1 / \lambda^{i}_{j}$ for four representative subdomains $\Omega_{j}$ is plotted in (a). The non-overlapping decomposition, with the four chosen subdomains highlighted, is presented in (b). In (c) and (d), a sample of eigenvectors $\phi^i_j$ for the subdomains $29$ and $56$ is illustrated, respectively. (For clarity, the first six eigenvalues for $\Omega_{j=55,56,64}^{\ast}$ are not plotted in (a), since by definition $\lambda^{0}_{j} = \ldots = \lambda^{5}_{j} = 0$.)}
    \label{fig:EValuesVectors}
\end{figure*}

We start by analysing the local generalized eigenvalue problems (GEVP), in particular the behavior of the smallest eigenvalue $\lambda^{m_j + 1}_{j}$ corresponding to any eigenvector not included in the local approximation space, which is the main parameter driving the method accuracy (cf.~the error bound in \cref{eq:error_bound}). First, the decay of $1 / \lambda^{i}_{j}$ for representative subdomains is analysed, as well as the shape of the associated eigenvectors. Then, the effect of the oversampling size $o^*$ is discussed. A particular focus will be on the local finite element aspect ratio which severely affects the accuracy of the computed eigenvectors/-values and thus also the observed actual decay of the error bound.

The beam domain in \Cref{fig:Models} (a) is decomposed into $64$ subdomains for the first experiment; see Fig.~\ref{fig:EValuesVectors} (b). In each subdomain, a local eigenvalue problem is solved to construct the local approximation space with a fixed number of $o^{\ast}=8$ oversampling layers. 
In \Cref{fig:EValuesVectors}~(c,d) a selection of exemplary, local \gls{GEVP} solutions on two representative subdomains are presented together with a semi-logarithmic plot of $1 / \lambda^{i}_{j}$ in Fig.~\ref{fig:EValuesVectors}(a). For the eigenvalue plot two further subdomains are added. In total, there is one subdomain intersecting $\Gamma_D$ ($j=29$), one intersecting $\partial \Omega$ $(j=64)$ and two interior subdomains ($j=55,56$), with the first one of the two having a higher surface-to-volume ratio.

The semi-logarithmic plot of $1 / \lambda^{i}_{j}$ in Fig.~\ref{fig:EValuesVectors}(a) demonstrates the predicted, nearly exponential decay of the local approximation error with respect to the basis size ($m_j$) in all cases.
The decay is faster for subdomains intersecting $\partial \Omega$; the more the subdomain intersects $\partial \Omega$ the faster is the decay. A further factor is the surface-to-volume ratio. The higher this ratio the slower the exponential decay of $1 / \lambda^{i}_{j}$, as exemplified by the relative decays of subdomains $j=29, 64, 55$ and $56$.
The partitioning of the domain could be optimised to unify the surface-to-volume ratio over all subdomains.
In fact, for simple, laminated composites a regular domain decomposition could be chosen, e.g., into rectangular subdomains. Here, however, we do not make this choice in order to show the robustness of the approach to rather general subdomain partitionings, as provided by automatic graph partitioners such as \emph{ParMetis}~\cite{karypis1998multilevelk}, and thus to show the potential of the approach for simulating very complex structures. 

In the interior subdomain $\Omega_{j=56}^{\ast}$, the first six eigenvectors (indices $0$ to $5$) depicted in \Cref{fig:EValuesVectors} (d) correspond to the zero energy modes (or rigid body modes) representing shifts and rotations of the structure. The following modes for $\Omega_{j=56}^{\ast}$ and the first few modes for $\Omega_{j=29}^{\ast}$ in \Cref{fig:EValuesVectors} (c) correspond to low-energy deformations of the subdomain: bending and shearing. The higher-energy modes correspond to higher frequency deformations ($i>50$).

\begin{figure}
    \centering
    \includegraphics[width=.7\linewidth]{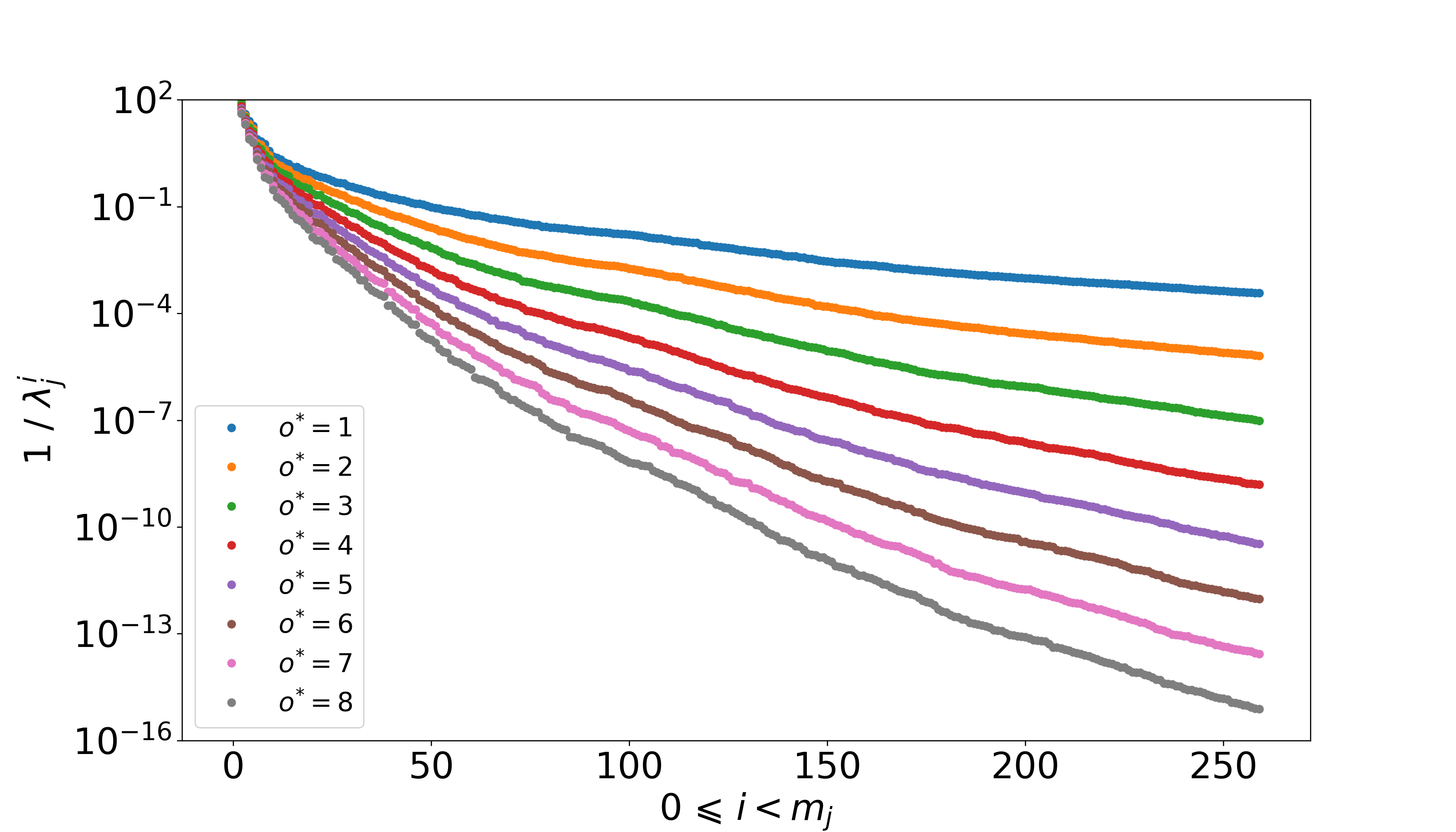}
    \caption{Effect of the amount of oversampling on the decay of the eigenvalues, and thus on the local error bound.}
    \label{fig:ovs_on_error_bound}
\end{figure}

As explained in \Cref{sec:geneo_as_gmsfem}, the oversampling parameter $o^{\ast}$ is key to the decay rate of the local approximation errors. In \Cref{fig:ovs_on_error_bound}, the reciprocal eigenvalues $1 / \lambda^{i}_{j}$ for subdomain $\Omega_{j=56}^{\ast}$ are presented 
for a range of oversampling sizes. The decay of $1 / \lambda^{i}_{j}$ clearly accelerates as the oversampling size is increased. Thus, for a desired error bound (i.e. $1 / \lambda^{i}_{j} < \num{e-7}$), subdomains with larger amounts of oversampling lead to significantly smaller local bases and consequently to an overall smaller coarse space. The trade-off, however, is that the dimensions of the local \gls{GEVP}s grow significantly with the amount of oversampling and this will be discussed later.
In general, the objective is to pick the minimal number of modes $m_j$ to approximate the solution (displacement) and its derivative (strain and stress) sufficiently well.
The accuracy of the coarse approximation and the effect of changing the main tuning parameters (basis / oversampling size) will be analysed in \Cref{sec:accuracy}.
The efficiency study will be performed for the more complex C-spar structure and is presented in \Cref{sec:scalability}.

\begin{figure}[t!]
    \centering
    \includegraphics[width=1.05\linewidth]{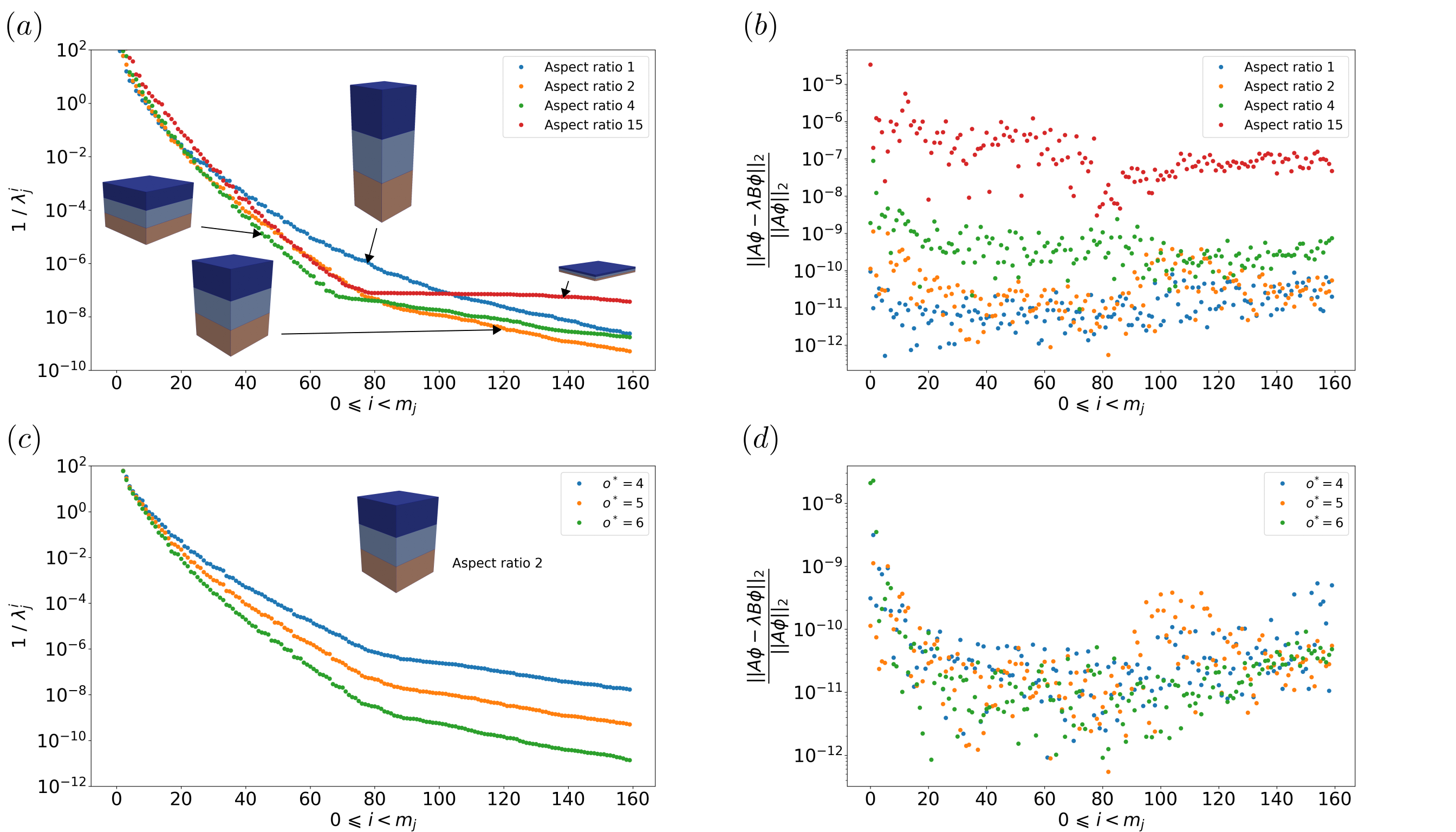}
    \caption{The reciprocals $1 / \lambda^{i}_{j}$ of the eigenvalues of the local GEVP (a,c) and the relative $L_2$-error in computing the associated eigenvectors (b,d); in particular, studying the effect of bad element aspect ratio (for $o^*=5$) in (a,b), as well as oversampling size for fixed aspect ratio in (c,d).}
    \label{fig:aspectratio}
\end{figure}

The results in Figures~\ref{fig:EValuesVectors} and \ref{fig:ovs_on_error_bound} are computed using square elements with an aspect ratio of one between horizontal (through-thickness) and vertical (in-plane) edges. However, such a choice leads to an unreasonably high number of elements in larger composite structures with a larger number of more realistic, thinner plies, such as the C-spar described in Section \ref{sec:compositeDescription} and depicted at he bottom of \Cref{fig:Models}.
To build a reasonably sized model, the in-plane discretization has to be reduced, leading to flat elements with a larger aspect ratio. Unfortunately this reduces the accuracy significantly, as presented for four element aspect ratios and $o^*=5$ in \Cref{fig:aspectratio} (b), where the relative $L_2$-error for each eigenpair is shown. As a consequence, after a similar initial decay (up to $i=60$) we observe a change in the slope of $1 / \lambda^{i}_{j}$ for aspect ratios bigger than one in \Cref{fig:aspectratio} (a). For very large aspect ratios of 15 and above, it even leads to a plateau in $1 / \lambda^{i}_{j}$ (see \Cref{fig:EValuesVectors} (a), red curve). This is to be expected, due to the larger condition numbers of the stiffness matrices in the GEVP leading to more unstable eigensolves.

Increasing the oversampling size does not alleviate this problem, as seen in \Cref{fig:aspectratio} (c,d); the relative $L_2$-errors in the eigenpairs are independent of $o^*$ and the slope of $1 / \lambda^{i}_{j}$ changes roughly at the same value of $i$. We also tested a different type of higher-order finite element, namely a 20-DoF quadratic serendipity element, but the loss of accuracy due to high element aspect ratios persists.
However, the slope change and the plateau of the eigenvalues for larger aspect ratios do not have a strong impact, since they occur only at values well below what is needed for good practial approximation, especially for higher oversampling sizes.

\subsection{Coarse approximation accuracy}\label{sec:accuracy}

In this section, we investigate the accuracy of the MS-GFEM approximation. 
The fine-scale reference solution is computed using an iterative \gls{CG} method with \gls{GenEO} as preconditioner. By setting a sufficiently low tolerance, the error due to the iterative solution via \gls{CG} can be neglected. 

We consider in the following the relative errors in $L_2$-norm between the coarse approximations of displacement and strain fields, $u_{H}$ and $\epsilon_{H}$, and their fine-scale counterparts, $u_{h}$ and $\epsilon_{h}$,
i.e.,
\[
e_d =  \dfrac{\lVert u_{h} - u_{H} \rVert_{2}}{\lVert u_{h} \rVert_{2}}, \qquad e_{\epsilon} =  \dfrac{\lVert \epsilon_{h} - \epsilon_{H} \rVert_{2}}{\lVert \epsilon_{h} \rVert_{2}}.
\]

\begin{figure}[t!]
    \centering
    \includegraphics[width=\linewidth]{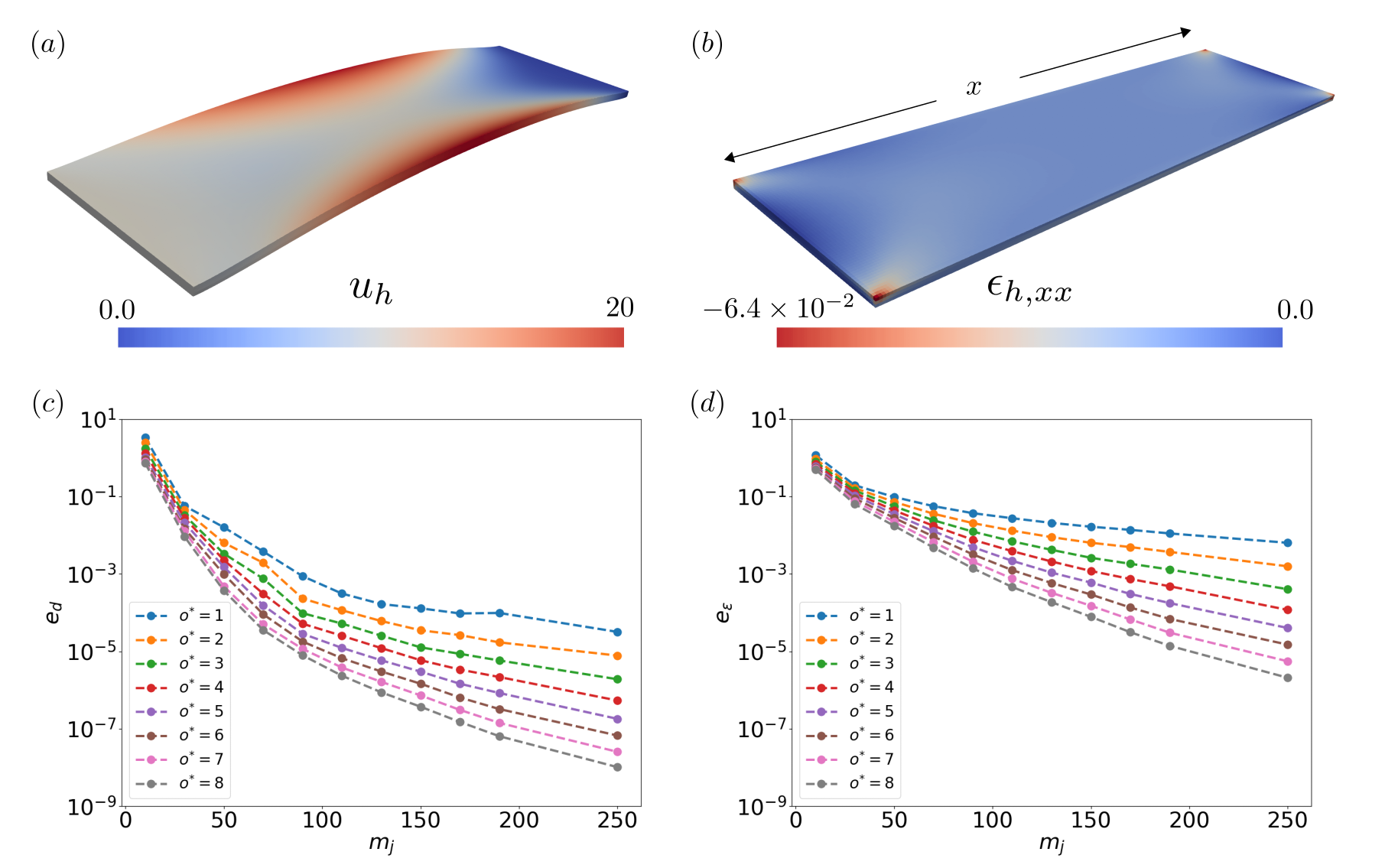}
    \caption{$L_2$-errors versus local basis size for the displacement field (c) and the strain field (d). The reference fine-scale approximations are displayed above in (a) and (b).}
    \label{fig:L2_error}
\end{figure}

\Cref{fig:L2_error} (a) shows the fine-scale approximation of the displacement field $u_{h}$ of the beam under compressive loading. The unit displacement (${\bf x}(z)=-10$ mm) applied here causes the characteristic out-of-plane deformation of the structure related to the non-symmetric stacking sequence chosen for this example.
A maximum displacement of $20$mm is observed on the two sides of the beam. \Cref{fig:L2_error} (b) depicts the fine-scale strain approximation in the direction of compression, denoted by $x$ here. The order of magnitude of strains is $\mathcal{O}(\num{e-2})$.

The main parameter driving the accuracy of the coarse approximation is the smallest eigenvalue $\min_j \lambda^{m_j+1}_{j}$ corresponding to any local eigenvector not included in the basis. As explained above, the size $m_j$ of each of the local bases and the oversampling size $o^{\ast}$ will be decisive to control this parameter. For our analysis we vary these two key parameters up to maxima of $m_j=250$ and $o^{\ast}=8$, for which the error bound lies below $\num{e-15}$.
The relative $L_2$-errors for displacement and strain are shown in \Cref{fig:L2_error}, resp.~(c) and (d), using a logarithmic scale. Both $e_d$ and $e_{\epsilon}$ decay exponentially with respect to $m_j$ for all amounts of oversampling. The decay rate of the errors is higher for larger oversampling sizes, which is in agreement with the behavior of the reciprocal eigenvalues depicted in \Cref{fig:ovs_on_error_bound}. 

\begin{figure}[t!]
    \centering
    \includegraphics[width=\linewidth]{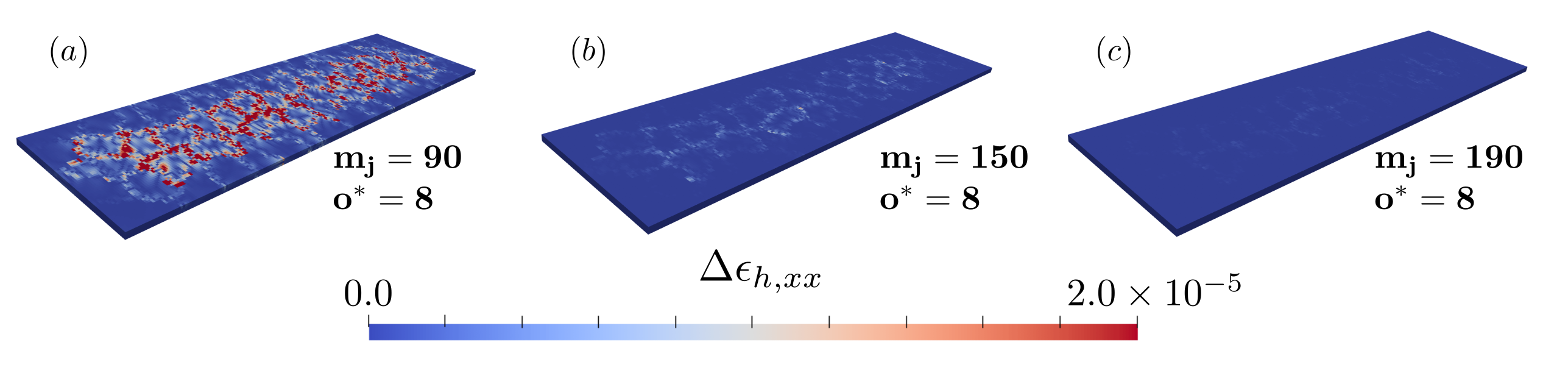}
    \caption{Localization of the error depending on the local space size $m_j$. The absolute difference $\Delta \epsilon_{h,xx}$ between the coarse approximation and the fine scale approximation is plotted.}
    \label{fig:local-error-mj}
\end{figure}
We further investigate the location of the largest strain error, in order to explain in more detail its dependence on $m_j$ and $o^{\ast}$. The coarse approximation of the strain field has been projected onto the fine-scale space to allow a direct comparison element-by-element. The absolute differences $\Delta \epsilon_{h,xx}$ between the fine-scale solution and the coarse approximation of the strain field are depicted in \Cref{fig:local-error-mj}
for three different values of the local basis size per subdomain, $m_j=\{90, 150, 190\}$, and a fixed value of $o^{\ast}=8$. The aim is to investigate the error behavior alongside the convergence curve for $o^{\ast}=8$ in \Cref{fig:L2_error} (d).
The three chosen basis sizes produce coarse-space solutions of high accuracy with relative $L_2$-errors bellow $\num{e-3}$. 
The error maxima are observed on subdomains located in the center of the beam. Indeed, the local error on each subdomain follows the decrease of $1 / \lambda_{j}^{i}$. As shown in \Cref{fig:EValuesVectors} (a), interior subdomains (e.g., $\Omega_{j=56}^{\ast}$) need more eigenvectors for the same local error. Eventually, a sufficiently large local basis size ensures an accurate approximation for all subdomains.

\begin{figure}[t!]
    \centering
    \includegraphics[width=\linewidth]{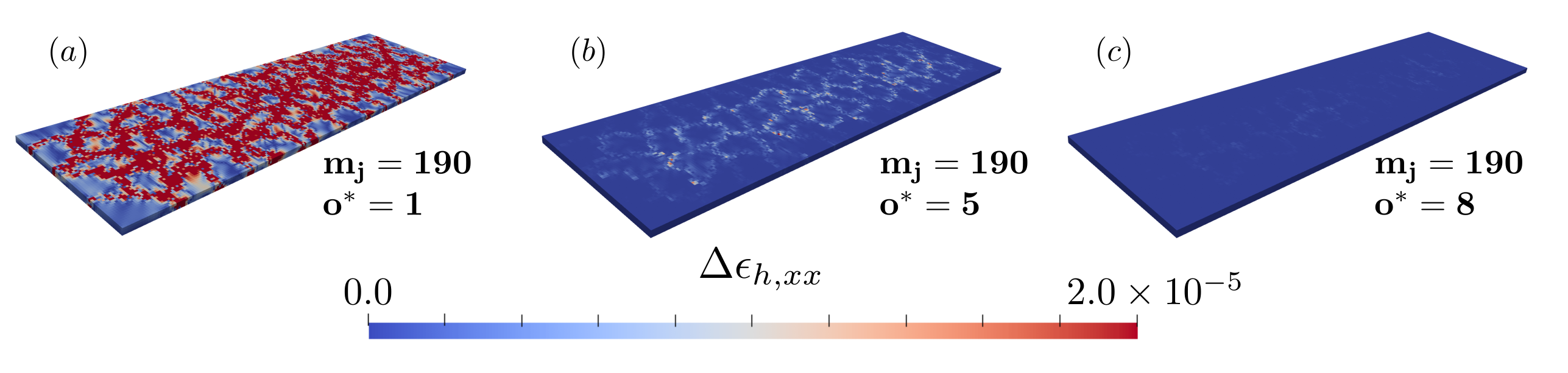}
    \caption{Localization of the error depending on the oversampling size $o^{\ast}$. The absolute difference $\Delta \epsilon_{h,xx}$ between the coarse approximation and the fine scale approximation is plotted.}
    \label{fig:local-error_ovs}
\end{figure}

Conversely, for fixed $m_j=190$ \Cref{fig:local-error_ovs} assesses the impact of the oversampling size, varying $o^{\ast} \in \{1, 5, 8\}$. As expected, small oversampling size leads to an error localised at subdomain interfaces.
Increasing the oversampling area decreases the interface errors, and, eventually a smooth and very low error is observed across the domain for $o^{\ast}>5$.
This example showcases one of the main advantages of the MS-GFEM: the interface problem at subdomain boundaries is handled by oversampling. Thus, the approximation is then robust to the chosen domain decomposition and there is no scale separation between the ply scale (mesoscopic) and the structural scale (macroscopic).

\begin{figure}
    \centering
    \includegraphics[width=\linewidth]{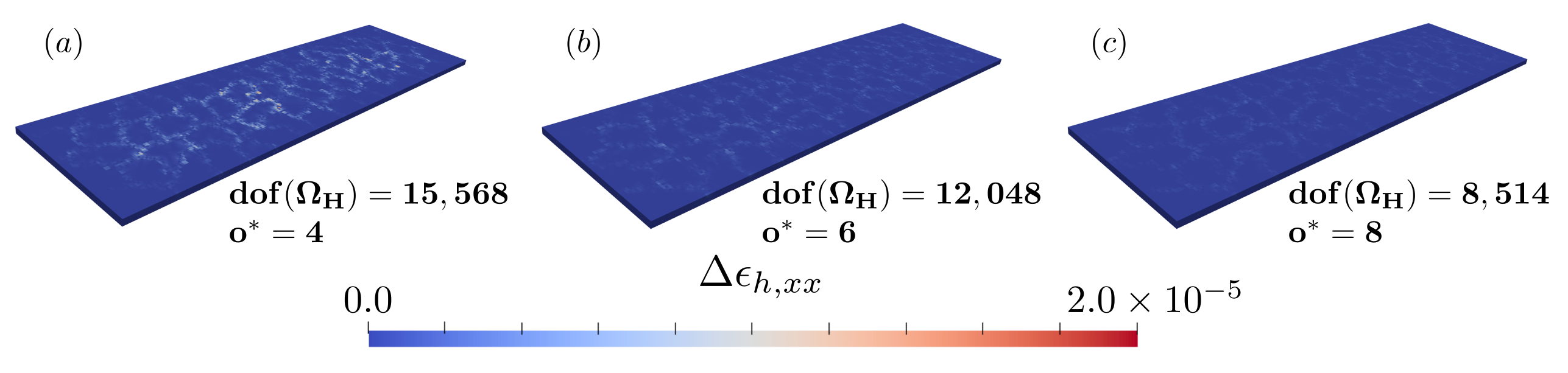}
    \caption{Localization of the error depending on the oversampling size $o^{\ast}$ for a given targeted accuracy. The absolute difference $\Delta \epsilon_{h,xx}$ between the coarse approximation and the fine scale approximation is plotted.}
    \label{fig:local-error_threshold}
\end{figure}

In \Cref{fig:local-error_threshold}, the strain error is evaluated for a fixed eigenvalue threshold $1 / \lambda_{j}^{m_j+1} = \num{e-7}$ for various oversampling sizes, to study the impact on the overall coarse space size.
For fixed accuracy, the size of the coarse space is reduced by $22.6\%$ when increasing $o^{\ast}=4$ to $o^{\ast}=6$ and by $45.3\%$ from $o^{\ast}=4$ to $o^{\ast}=8$. This improved model order reduction obviously reduces the cost of the coarse problem solve, but it does increase the size of the local problems and thus the cost of solving the local GEVPs. This cost trade-off will be studied more carefully in \Cref{sec:scalability} for the C-spar.

\subsection{Importance of the A-harmonic condition for composite problems}\label{sec:comparison_two_coarse spaces}

\begin{figure}[t!]
    \centering
    \includegraphics[width=\linewidth]{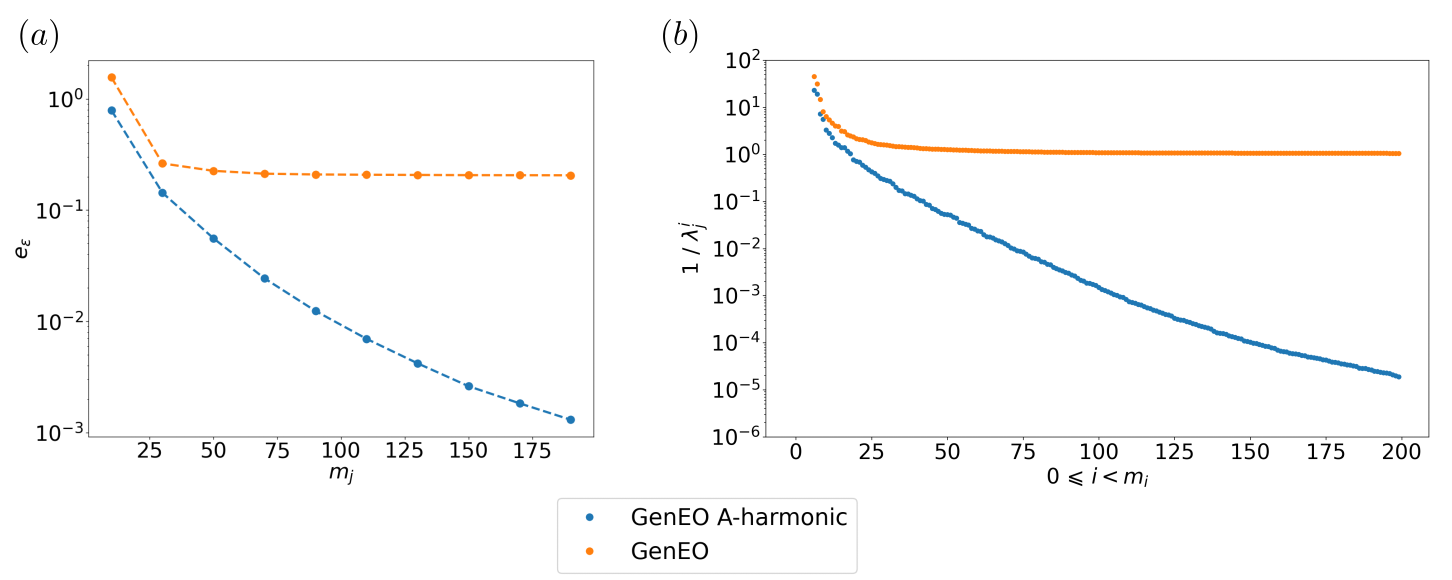}
    \caption{Accuracy of the GenEO-type coarse spaces. Absolute relative $L_2$-error on the strain field for both the classical GenEO and the A-harmonic GenEO (MS-GFEM) approximations (a). The quotient $1 / \lambda^{i}_{j}$ for one particular subdomain for each of the approaches (b).}
    \label{fig:error-geneo}
\end{figure}

To emphasize the necessity of the A-harmonic condition in our framework, we compare the new MS-GFEM space to the classic GenEO coarse space with the GEVPs formulated directly in the \gls{FE} space $V_h$, i.e., without enforcing A-harmonicity.
The theoretical error bound in Theorem \ref{thm:error_bound} 
is the same for both GFEM spaces, with and without the A-harmonic condition in the GEVP. Crucially, in both cases the approximation error is bounded by the reciprocal $1 / \lambda^{i}_{j}$ of the smallest eigenvalue corresponding to any eigenvector that is not included in the coarse space $V_H$.

\Cref{fig:error-geneo} (right) shows that the decay of $1 / \lambda^{i}_{j}$ for $i<10$ is comparable between classic GenEO and the A-harmonic formulation, since the first modes (rigid body, shear and bending) appear in both. In fact, the classic GenEO coarse space with a small number of lowest-eigenvalue modes was shown in Reinarz et al.~\cite{butler2020high} to provide a robust preconditioner within CG that reduces the condition number effectively and leads to a low number of \gls{CG} iterations for composites problems.
Beyond that, we observe in \Cref{fig:error-geneo} (right) that the spectrum for the classic GenEO formulation decays much slower and eventually almost stagnates, when compared to MS-GFEM. This is reflected in the coarse approximation error in \Cref{fig:error-geneo} (left), as predicted in Theorem \ref{thm:error_bound}.
Even including $190$ eigenvectors in classic GenEO is not sufficient to accurately represent the strain field in this application (with $e_{\epsilon} \approx 1$). MS-GFEM with A-harmonic GEVPs, on the other hand, achieves an excellent coarse approximation at much lower basis sizes.

\subsection{Method accuracy on an aerospace part}\label{sec:csparaccuracy}

\begin{figure}[t!]
    \centering
    \includegraphics[width=0.9\linewidth]{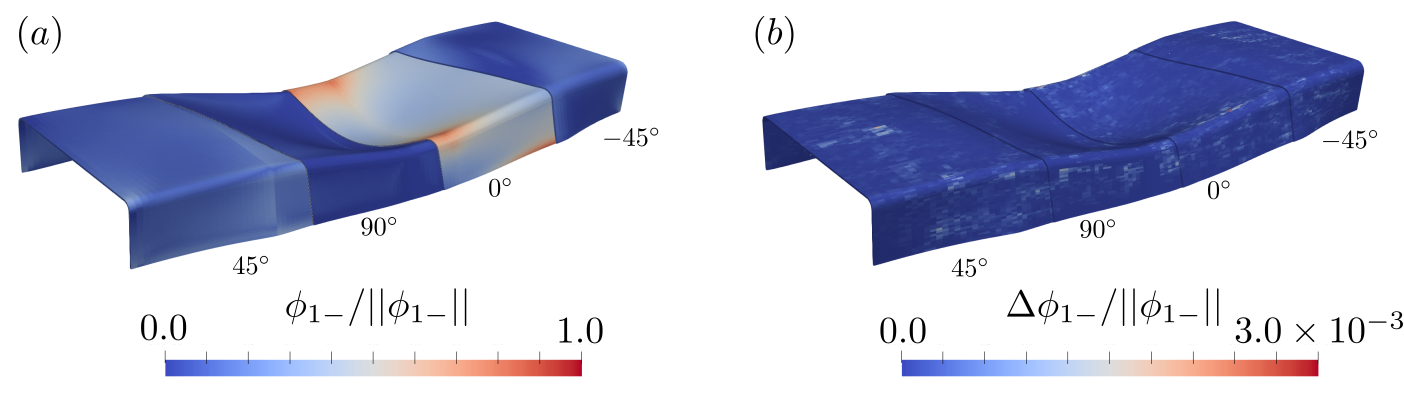}
    \caption{Application of the method to the C-spar. The elementwise, normalized compressive failure criterion is shown in (a) and the error $\Delta \phi_{1-}$ of the MS-GFEM approximation of the failure criterion with respect to a fine scale computation of the failure criterion in (b).}
    \label{fig:cspar-application}
\end{figure}
In this section, the behavior of the C-spar, described in \Cref{sec:compositeDescription}, under compression is investigated. The domain, shown in \Cref{fig:Models}~(b), is divided into $256$ subdomains. For an 
oversampling size of $o^{\ast}=4$ and a tolerance of $t = \num{e-6}$ on $1/\lambda_j^{m_{j}+1}$ we obtain local basis sizes $m_j$ varying between $60$ and $170$. These parameters lead to an average subdomain size of $\text{dof}(\Omega_{j}^{\ast}) = 30000$ and a model order reduction of $32\times$. With 4GB of RAM per subdomain to compute the local GEVP this parameter choice requires 4 nodes of the HPC cluster \textit{Hamilton8} for the 256 subdomains. The element aspect ratio is 15 for this example (see \Cref{fig:Models}). Thus, as for the simple beam example, the possible coarse approximation accuracy is limited by the fine scale error.

In order to assess the coarse approximation, the longitudinal compressive failure criterion
\begin{equation}\label{eq:criterion}
    \phi_{1-} = \frac{ \left\langle |\tilde{\sigma}^{R}_{12}| +\eta^L  \tilde{\sigma}^{R}_{22} \right\rangle }{S_L}
\end{equation}
from~\cite{furtado2019simulation} is computed, which is a linear combination of relevant stress components: transverse and shear. A detailed description of the material parameters that are involved is outlined in \cite{furtado2019simulation}; the notation $\langle \cdot \rangle$ denotes the positive part of the argument. 
It represents the initiation of the fiber kinking mechanism, observed for fiber reinforced composites under compression.
In slender composite structures, this failure mode typically appears after the buckling of the structure, but even in the pre-buckling context, this criterion gives a good insight into the usefulness of the method for studying compressive failure.

The elementwise longitudinal compressive failure criterion, computed in local coordinates, is plotted in \Cref{fig:cspar-application}~(a). To compute it the coarse approximation $u_H$ is first projected onto the fine scale FE space $V_h$. The relative error with respect to a direct meso-scale approximation in $V_h$ is shown in \Cref{fig:cspar-application} (b).
With the chosen parameters, the relative error on the failure criterion is below $\num{3.0e-3}$. 
In \Cref{fig:cspar-application} (a), four plies (ply 6, 12, 18 and 24) are highlighted in a certain part of the domain, representing each stacking orientation and showing that the coarse approximation is able to accurately represent the quick variation of the compressive failure criterion through thickness.
The remaining error is small and uniformly distributed, with higher values near subdomain boundaries and on interior subdomains, in agreement with the observations on the beam example. 

\subsection{Method scalability}\label{sec:scalability}

\begin{figure}[t!]
    \centering
    \includegraphics[width=\linewidth]{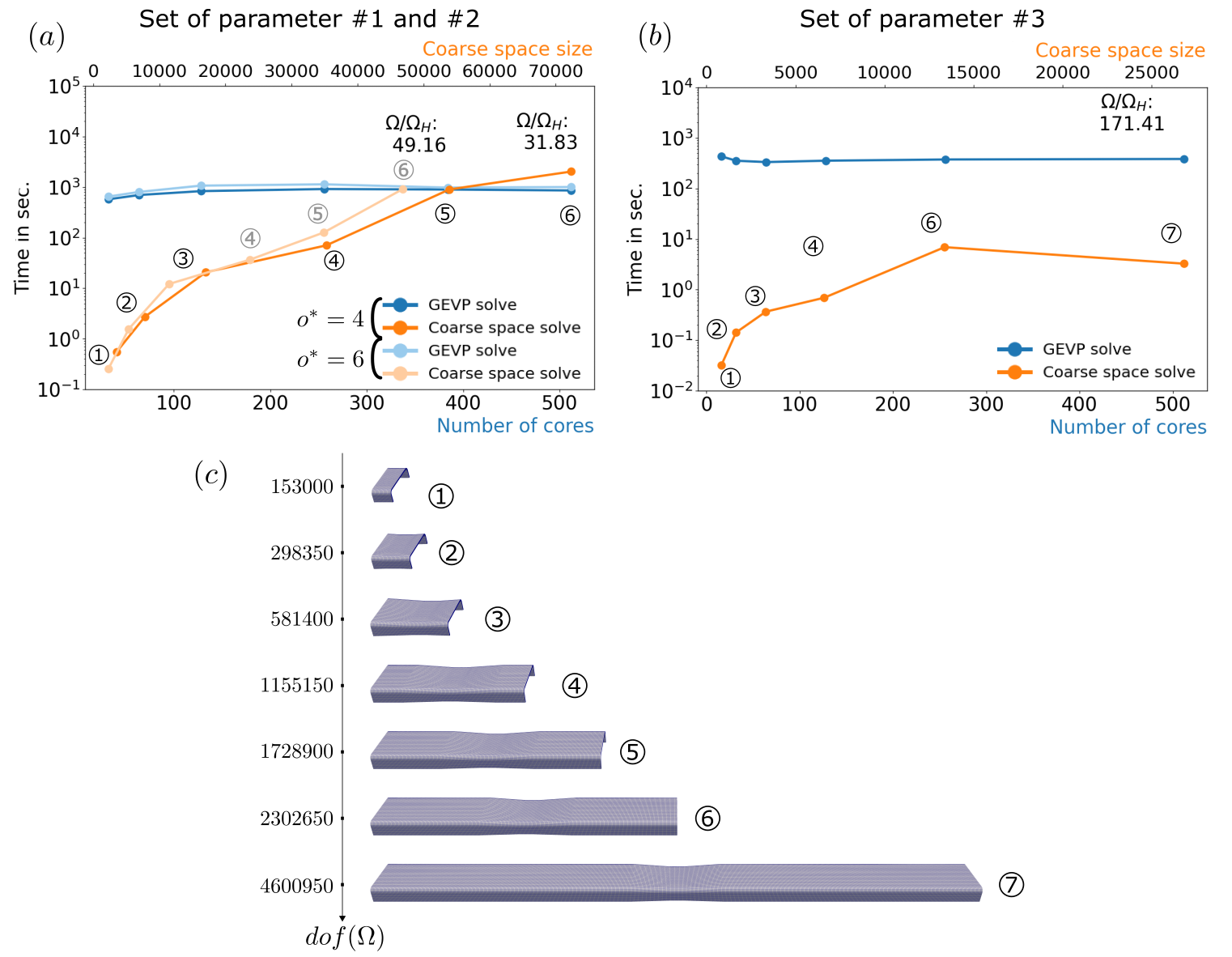}
    \caption{Scaling of the method applied to C-spars of different lengths for three sets of parameters. In (a) the Sets $\#1$ and $\#2$ optimize the approximation accuracy, while in (b) Set $\#3$ aims at a very efficient simulation with acceptable accuracy. The FE grids used are displayed in (c).}
    \label{fig:scalability}
\end{figure}

A parallel (weak) scaling test of the method is presented in \Cref{fig:scalability}, where the cost of the local GEVP solves and the cost of the coarse solve are assessed using C-spars of different lengths while fixing the discretization of the meso-structure: $24$ elements through thickness (one per ply) and a constant element size in the other directions (aspect ratio 15). The C-spar models vary between $L=62.5$mm with $\text{dof}(V_h)=\num{1.53e5}$ (Model 1) and $L=2$m with $\text{dof}(V_h)=\num{4.6e6}$ (Model 7). The number of processors ($P$) employed for each model have been chosen such that $\text{dof}(V_h)/P$ remains constant.

Three combinations of parameters have been tested, choosing $o^*=4$ and $o^*=6$ and selecting a threshold of $t=\num{e-6}$ and $t=\num{e-3}$ for the reciprocal of the local eigenvalues. Parameter sets $\#1$ and $\#2$  assess the scalability of the method for a highly accurate solution with $t=\num{e-6}$. In Set $\#1$ $o^{\ast}=4$, chosen to reduce RAM consumption, whereas in Set $\#2$ a larger oversampling size of $o^{\ast}=6$ is selected, thus requiring more RAM. This trade-off will be further discussed below. Set $\#3$ has been selected to study how the method scales for a lower eigenvalue threshold of $t=\num{e-3}$ with smaller $o^{\ast}=3$ and half the amount of subdomains in Sets $\#1$ and $\#2$, such that $\text{dof}(\Omega_{j}^{\ast})$ is roughly the same. 

For all three sets, the GEVP solve step scales perfectly, since it requires no parallel communication. As expected the local GEVP solve times increase slightly between Sets $\#1$ and $\#2$ due to the slightly bigger local problem sizes. The overall cost for the coarse space setup for each model is essentially identical, since it is dominated by the cost of solving the local GEVPs and can also be carried out fully in parallel with only some local data exchange. In contrast, the final coarse space problem \eqref{eq:GFEM_coarse_approximation} is (currently) solved on one processor, and thus eventually dominates the overall cost for larger models. To remain efficient, the coarse space size needs to be optimized. For all considered models, the model order reduction is around $32\times$ for Set $\#1$ and around $50\times$ for Set $\#2$. As a consequence, the increase in oversampling size from Set $\#1$ to Set $\#2$ significantly reduces the cost of the coarse solve, particularly for the larger models (4–5–6). 
These gains clearly outweigh the higher cost for the local GEVP solves, but there is a limit to this improvement, preventing the use of too large oversampling sizes, namely the amount of local memory (RAM) needed to process the local GEVPs in parallel. 
In this example, subdomains with $o^{\ast}=6$ have more than $50,000$ DoF, necessitating over 8GB of RAM per subdomain for the local GEVPs. A reasonable oversampling size thus needs to balance performance and memory consumption.

\begin{figure}[!t]
    \centering
    \includegraphics[width=.8\linewidth]{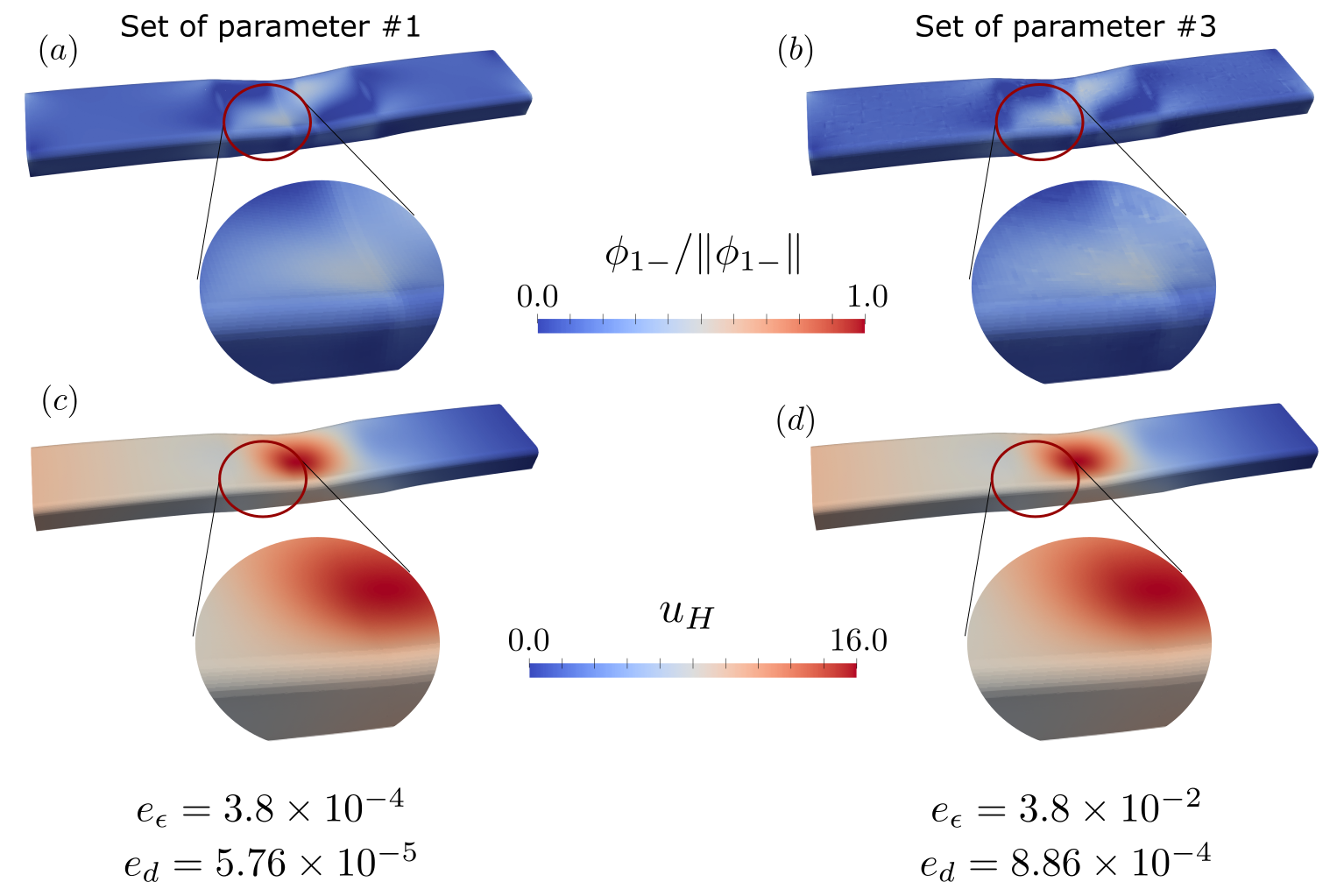}
    \caption{Visualisation of the solution accuracy for
    Set $\#3$ (b,d) compared to Set $\#1$ (a,c) for Model 6 (see \Cref{fig:scalability}), plotting for both sets in (a, b) the compressive criterion (\cref{eq:criterion}) and in (c, d) the displacement field. The relative $L_2$-errors (on strain  $e_{\epsilon}$ and displacement $e_{d}$) with respect to the fine scale approximation are associated with each parameter set.}
    \label{fig:scalabilityp2}
\end{figure}

The strength of our approach is the ability to construct approximation spaces of adjustable complexity in a very simple way. A cheaper approximation can be built by choosing a lower threshold $t$ for selecting the local eigenvectors. In particular, for Set $\#3$, where $o^{\ast}=3$ and $t=\num{e-3}$, the size of local bases $(m_{j})$ reduces to between $30$ and $50$, leading to a 
model order reduction of about $170\times$. This allows the construction and solve of the coarse-space approximation for a $2$m C-spar (Model 7) on 512 processors in 10 minutes, further pushing the scalability test. 
The model order reduction is then sufficient to keep all the coarse solve computation times under $10$ seconds, even for the Model 7. 
The accuracy of the solution of two sets of parameters: $\#1$ and $\#3$ is compared in \Cref{fig:scalabilityp2} for Model 6.
Both qualitatively, as visible on \Cref{fig:scalabilityp2}, and quantitatively (relative $L_2$-errors) the displacement field $u_H$ appears to be well represented for both sets, which is in agreement with the beam example observations. 

However, as expected, since it is derived from the stresses (\cref{eq:criterion}), the accuracy of the compressive criterion in Set $\#3$ is affected by the small basis size.
Despite a visible noise on the compressive failure criterion of Set $\#3$ (see the zoom-in in \Cref{fig:scalabilityp2}), global extrema are preserved. Hence, this set of parameters is able to detect the global maximum of the criterion at very cheap cost.

\section{Conclusion \& future work}
In this paper, we have presented the first scalable \gls{HPC} implementation of a \gls{MS-GFEM} method and demonstrated that it delivers high quality approximate solutions for very small coarse space sizes. As proven in previous theoretical work, this is due to the nearly exponential decay of the reciprocal eigenvalues in the local generalised eigenvalue problems. Here, we have demonstrated that this nearly exponential decay crucially relies on enforcing an A-harmonic constraint on the local eigenproblems also in composites applications, and that oversampling of the local subdomains is essential to achieve good accuracies at small local basis sizes.
While the related GenEO-coarse space, which does not enforce A-harmonicity in the local eigenproblems, also leads to acceptable results in approximating displacements, we have seen that A-harmonicity is crucial to accurately approximate strains, stresses and derived failure criteria in composite applications.

We have demonstrated good parallel scalability on several hundreds of processor cores. While a single solve of the fine-scale problem is cheaper using, e.g., \gls{GenEO}-preconditioned Krylov methods, the localized approach of MS-GFEM opens up new opportunities for parallel scalability. When solving large numbers of closely related problems, eigenvectors from unaffected subdomains may be retained, solving costly eigenproblems only where model parameters or geometry changes between runs. If only a few subdomains are affected, the global solution can be computed using a significantly smaller number of processors in the same time as a full run. This is especially interesting in future Uncertainty Quantification (UQ) applications, such as the impact of (meso-scale) localized wrinkles in composite structures on the strength or the failure behaviour of the overall (macro-scale) structure. In such applications, the problem setup will essentially be identical in all but a few subdomains and large numbers of runs are required. 

The integration of our new \gls{MS-GFEM} method into an offline/online framework, where local approximation spaces will only be updated in a few subdomains between runs, is currently ongoing and will form part of a subsequent publication. This will also include the application of the offline/online framework as part of Uncertainty Quantification (UQ) methods for composites; in particular, exploiting the natural hierarchy of approximate models in the MS-GFEM framework
within multilevel UQ methods such as MLMC~\cite{MLMC} or MLMCMC~\cite{dodwell2019multilevel}. 

The coarse space solves have been handled by a single processor in this study. Another way to improve the framework efficiency will be to parallelize the coarse space solve, using a direct or an iterative parallel solver.
This aspect and the management of parallel resources in the online runs will also be explored.

On the practical side, our approach has a largely automatic workflow, from domain decomposition to the automatic generation of a coarse space specifically tuned to the given problem. The balance between global approximation error and basis size can be controlled by setting a single threshold for the selection of eigenvectors.
In multiscale applications, and specifically composites, MS-GFEM is particularly interesting since we obtain a low-dimensional approximation without assuming scale separation. Instead, the eigenproblems capture the structure of the given problem, providing better quality than hand-tuned approximations. The resulting coarse space then accurately captures fine- and coarse-scale interaction.

A very relevant aspect in the study of composite materials is material failure under load. In contrast to linear elasticity (as covered by this work), non-linear models are needed to simulate the failure of composites aero-structures (non-linear geometry, damage initiation and propagation). We are therefore also extending our methods to nonlinear solvers and implement nonlinear material behavior in \texttt{dune-composites}.

\section{Acknowledgement}
\label{sec:acknowledge}

The research was supported by the UK Engineering and Physical Sciences Research Council (EPSRC) through the Programme Grant: ‘Certification of Design: Reshaping the Testing Pyramid' EP/S017038/1 (https://www.composites-certest.com/). This multidisciplinary project aims at developing new approaches to enable the design and certification of lighter, more cost and fuel efficient composite aero-structures. The funding received is gratefully acknowledged. This work made use of the facilities of the Hamilton HPC Service of Durham University.
We thank Anne Reinarz (Durham University) for her support in the HPC experiments.
Richard Butler holds a Royal Academy of Engineering-GKN Aerospace Research Chair in Composites.



 \bibliographystyle{elsarticle-num} 
 \bibliography{cas-refs}





\end{document}